\documentclass[english]{amsart}
\usepackage[bookmarksnumbered,plainpages,hypertex]{hyperref}
\usepackage{amsfonts}
\usepackage{graphicx}
\usepackage[all]{xy}
\usepackage{srcltx}

\allowdisplaybreaks
\usepackage{babel}

\numberwithin{equation}{section} 
\numberwithin{figure}{section} 
\newtheorem{thm}{Theorem}[section]
\newtheorem{prop}[thm]{Proposition}
\newtheorem{lem}[thm]{Lemma}
\newtheorem{cor}[thm]{Corollary}
\theoremstyle{definition}
\newtheorem{defn}[thm]{Definition}
\newtheorem{example}[thm]{Example}
\newtheorem{claim}[thm]{}
\newtheorem{rem}[thm]{Remark}

\begin{document}

\title{A Milnor-Moore Type Theorem for Primitively Generated Braided Bialgebras}

\author{Alessandro Ardizzoni}

\curraddr{University of Ferrara, Department of Mathematics, Via Machiavelli
35, Ferrara, I-44121, Italy}

\email{alessandro.ardizzoni@unife.it}

\urladdr{http://www.unife.it/utenti/alessandro.ardizzoni}


\thanks{This paper was written while the author was member of GNSAGA with
partial financial support from MIUR within the National Research Project
PRIN 2007.}

\begin{abstract}
A braided bialgebra is called primitively generated if it is generated as an algebra by its space of primitive elements. We prove that any primitively generated braided bialgebra is isomorphic to the universal enveloping algebra of its infinitesimal braided Lie algebra, notions hereby introduced. This result can be regarded as a Milnor-Moore type theorem for primitively generated braided bialgebras and leads to the introduction of a concept of braided Lie algebra for an arbitrary braided vector space.
\end{abstract}


\maketitle
\tableofcontents

\section{Introduction}

\textbf{The main Idea.} Le $K$ be a fixed field. All vector spaces are understood as vector spaces over
$K$. Let $H$ be a pointed Hopf algebra (i.e. all its simple
subcoalgebras are one-dimensional) and denote by $G$ the set of grouplike
elements in $H$. Let $\mathrm{gr}H$ be the graded coalgebra associated
to the coradical filtration of $H$. It is well known that $\mathrm{gr}H$ is
isomorphic, as a Hopf algebra, to the Radford-Majid bosonization $R\#KG$
of $R$ by the group algebra $KG$, where $R$, the so-called diagram
of $H$, is a suitable connected braided bialgebra in the braided monoidal
category of Yetter-Drinfeld modules over $KG$. This is the starting
point of the \emph{lifting method} for the classification
of finite dimensional pointed Hopf algebras, introduced by N. Andruskiewitsch
and H.-J. Schneider, see e.g. \cite{AS- Lifting}. Accordingly to this method, first one has to describe
$R$ by generators and relations, then to lift the informations obtained
to $H.$ It is worth to notice that Andruskiewitsch and Schneider conjectured that, in the finite dimensional case and characteristic zero, $R$ is always primitively generated, \cite[Conjecture
1.4]{AS- FiniteQuantCartan} and \cite[Conjecture 2.7]{AS}. 

Therefore, in many cases, for proving certain properties of Hopf algebras,
one can reduce to the connected case. The price to pay is to work with bialgebras in a braided category, and not with ordinary bialgebras. 

Actually in our case it is more convenient to
work with braided bialgebras, that were introduced in \cite{Ta}.
Recall that a \emph{braided vector space} $(V,c)$ consists of a vector
space $V$ and $K$-linear map $c:V\otimes V\rightarrow V\otimes V$,
called braiding, satisfying $c_{1}c_{2}c_{1}=c_{2}c_{1}c_{2}$ (the
so-called quantum Yang-Baxter equation). Here $c_{1}=c\otimes V$
and $c_{2}=V\otimes c.$ A \emph{braided bialgebra} is then a braided
vector space which is both an algebra and a coalgebra with structures
suitably compatible with the braiding. Examples of braided bialgebras
are all bialgebras in those braided monoidal categories which are monoidal subcategories of the category of vector spaces.

A coalgebra (a fortiori a braided bialgebra) is called connected if
its coradical (i.e. the sum of all its simple subcoalgebras) is one-dimensional.
The main result of this paper, Theorem \ref{teo: generated}, states
that every primitively generated (whence connected) braided bialgebra
is isomorphic, as a braided bialgebra, to the universal enveloping algebra
of its infinitesimal braided Lie algebra. This result can be seen
as an extension of the celebrated Milnor-Moore Theorem \cite[Theorem 5.18]{Milnor-Moore}
for cocommutative connected bialgebras (once observed that such a
bialgebra is always primitively generated): in characteristic zero,
any cocommutative connected bialgebra is the enveloping algebra of its
space of primitive elements, regarded as a Lie algebra in a canonical way. 

In order to understand the statement of the main theorem, we need
to clarify what is meant by infinitesimal braided Lie algebra
of a braided bialgebra and by universal enveloping
algebra of it. The definition of infinitesimal braided Lie
algebra requires that of the braided tensor algebra $T(V,c)$ associated
to a braided vector space $\left(V,c\right)$. As an algebra, $T(V,c)$ 
is the tensor algebra $T=T\left(V\right)$. Then $T$ can be endowed
with a bialgebra structure depending on $c$ which is denoted by $T(V,c)$
and called the braided tensor algebra of $(V,c)$. 

Now, let $A$ be a braided bialgebra. The braiding of $A$ induces a braiding
$c_{P}$ on the space $P=P\left(A\right)$ of primitive elements in
$A$. The braided vector space $\left(P,c_{P}\right)$ will be called
the infinitesimal part of $A$.

Let us construct a direct system of braided bialgebras \[
U^{\left[0\right]}\overset{\pi_{0}^{1}}{\rightarrow}U^{\left[1\right]}\overset{\pi_{1}^{2}}{\rightarrow}U^{\left[2\right]}\overset{\pi_{2}^{3}}{\rightarrow}\cdots\]
The first term is $U^{\left[0\right]}:=T\left(P,c_{P}\right).$ The
usual universal property of $T(P)$ yields an algebra map $\varphi^{\left[0\right]}:U^{\left[0\right]}\rightarrow A$,
which turns out to be a braided bialgebra map, such that the diagram \[
\xymatrix@R=10pt{P\ar[dr]_{i}\ar[rr]^{i_{U^{[0]}}} &  & U^{[0]}\ar@{.>}[dl]^{\varphi^{[0]}}\\
 & A}
\]
commutes, where $i_{U^{[0]}}$ is the canonical map and $i$ is the inclusion. Set $P^{[0]}:=P\left(U^{\left[0\right]}\right)$. Then $i_{U^{[0]}}$
corestricts to $i^{[0]}:P\rightarrow P^{[0]}$. Moreover, $\varphi^{\left[0\right]}$ preserves primitive
elements so that it induces a map $b^{\left[0\right]}:P^{[0]}\rightarrow P$.
The second term $U^{\left[1\right]}$ of the direct system is then defined as the quotient of $U^{[0]}$ by its two-sided ideal generated by $\left[\mathrm{Id}-i^{\left[0\right]}b^{\left[0\right]}\right]\left[P^{[0]}\right]$. From the definition one gets $\varphi^{[0]}i^{[0]}b^{[0]}=\varphi^{[0]}$
so that $\varphi^{\left[0\right]}$ quotients to a bialgebra homomorphism
$\varphi^{\left[1\right]}:U^{\left[1\right]}\rightarrow A$ such that the diagram
\[
\xymatrix@R=10pt{U^{[0]}\ar[dr]_{\varphi^{[0]}}\ar[rr]^{\pi_{0}^{1}} &  & U^{[1]}\ar@{.>}[dl]^{\varphi^{[1]}}\\
 & A}
\]
commutes, where $\pi_{0}^{1}$ denotes the canonical projection. Set $P^{[1]}:=P\left(U^{\left[1\right]}\right)$.
Then $i_{U^{\left[1\right]}}=\pi_{0}^{1}i_{U^{\left[0\right]}}$ corestricts
to $i^{[1]}:P\rightarrow P^{[1]}$. Moreover, $\varphi^{\left[1\right]}$
induces a map $b^{\left[1\right]}:P^{[1]}\rightarrow P.$ The third term $U^{\left[2\right]}$ of the direct system is then the quotient of $U^{[1]}$ by its two-sided ideal generated by $\left(\mathrm{Id}-i^{\left[1\right]}b^{\left[1\right]}\right)\left[P^{[1]}\right]$. Moreover $\varphi^{\left[1\right]}$ quotients to a bialgebra homomorphism $\varphi^{\left[2\right]}$ such that $\varphi^{\left[2\right]}\pi_1^2=\varphi^{\left[1\right]}$.

Going on this way, we obtain the direct
system $\left(U^{\left[n\right]}\right)_{n\in\mathbb{N}}$
and a sequence of maps $b_{P}:=\left(b^{\left[n\right]}\right)_{n\in\mathbb{N}}$, where $b^{\left[n\right]}:P^{[n]}\rightarrow P.$ 
The datum $(P,c_{P},b_{P})$ is called the infinitesimal braided Lie
algebra of $A$. 

The direct limit of the direct system $\left(U^{\left[n\right]}\right)_{n\in\mathbb{N}}$
is a braided bialgebra that we denote by \[
U\left(P,c_{P},b_{P}\right):=\underrightarrow{\lim}U^{\left[n\right]}.\]
The compatible family of morphisms $\left(\varphi^{\left[n\right]}\right)_{n\in\mathbb{N}}$
gives rise to a braided bialgebra homomorphism $\varphi^{\left[\infty\right]}:U\left(P,c_{P},b_{P}\right)\rightarrow A$.
We can show that the space of primitive elements in $U\left(P,c_{P},b_{P}\right)$ identifies with $P$. Using this fact, in Theorem \ref{teo: generated}, by applying a famous result due to Heyneman and Radford (see \cite[Theorem 5.3.1]{Montgomery}), we can prove that $\varphi^{\left[\infty\right]}$
is an isomorphism whenever $A$ is primitively generated. Now, the
structure and properties of $U\left(P,c_{P},b_{P}\right)$ are encoded
in the datum $\left(P,c_{P},b_{P}\right).$ Indeed this leads to the
introduction of what will be called a braided Lie algebra $\left(V,c,b\right)$
for any braided vector space $\left(V,c\right)$ and of the related
universal enveloping algebra $U\left(V,c,b\right)$.\medskip

\textbf{The results in Detail.} The paper is organized as follows. 

Section \ref{sec:prelim} contains preliminary facts and notations that will be used in the paper. 

Section \ref{sec:pre-categ} concerns general results on pre-categorical and categorical spaces that will be applied in next section. 

Section \ref{sec:brideal} deals with the problem of quotienting braided bialgebras. In particular, we investigate braided bialgebra quotients of the braided tensor algebra $T(V,c)$ associated to a given braided vector space $(V,c)$. In Lemma \ref{lem:quotOfTensor}, these braided bialgebras are characterized as $(V,c)$-bialgebras, a notion hereby introduced. In Theorem \ref{teo: (V,c)-bialg}, we show how to construct a quotient of a $(V,c)$-bialgebra whenever there exists a suitable map. This fact will be relevant in defining a braided Lie algebra of a braided vector space.

Section \ref{sec:univenv}, is devoted to the introduction and investigation of the notion of braided Lie algebra and universal enveloping algebra. Let $\left(V,c\right)$ be a braided vector space. We say that
$\left(V,c\right)$ has a bracket if there is a family of maps $b=\left(b^{\left[n\right]}:P^{[n]}\rightarrow V\right)_{n\in\mathbb{N}}$
which allow, roughly speaking, to construct step by step a direct
system of braided bialgebras \[
T(V,c)=U^{\left[0\right]}\overset{\pi_{0}^{1}}{\rightarrow}U^{\left[1\right]}\overset{\pi_{1}^{2}}{\rightarrow}U^{\left[2\right]}\overset{\pi_{2}^{3}}{\rightarrow}\cdots\]
as we did in the case of $(P,c_{P})$. The direct limit of this direct
system is a braided bialgebra denoted by \[
U\left(V,c,b\right):=\underrightarrow{\lim}U^{\left[n\right]}.\]
This bialgebra is called the universal enveloping algebra of $(V,c,b).$
We say that $\left(V,c,b\right)$ is a \emph{braided Lie algebra} whenever
$\left(V,c\right)$ is a braided vector space, $b$ is a bracket on
$\left(V,c\right)$ and the canonical map $i_{U}:V\rightarrow U\left(V,c,b\right)$
is injective (the latter will be called the \emph{implicit
Jacobi identity}). 
In Corollary \ref{cor:P(U)}, we prove that, for a braided Lie algebra $(V,c,b)$, the space of primitive elements in $U(V,c,b)$ identifies with $V$ through the canonical map $i_U:V\rightarrow U$.
In Proposition \ref{pro: i_U split}, we get that for a braided vector space $(V,c)$ with a bracket $b$, the implicit Jacobi identity holds if and only if $b^{[n]}i^{[n]}=\mathrm{Id}$ for all $n\in\mathbb{N}$. Applying this result, in Corollary \ref{coro: ker b}, we give the alternative description of $U^{[n+1]}$ as the quotient of $U^{[n]}$ by its two-sided ideal generated by $\mathrm{Ker}(b^{[n]})$.

Let $(V,c,b)$ be a braided Lie algebra and let $U=U(V,c,b)$. In Theorem \ref{teo:Ardi}, we associates to $(V,c,b)$ a braided Lie algebra $(V,c,\beta)$ in the sense of \cite{Ardizzoni-Universal}. Furthermore, we show that the corresponding universal enveloping algebra $\mathbb{U}(V,c,\beta)$, also defined in \cite{Ardizzoni-Universal}, coincide with $U^{[1]}$ as given above. Therefore, the notions defined in \cite{Ardizzoni-Universal} can be seen as a first step of the procedure we consider in the present paper. Now, for a given braided vector space $(V,c)$, let $E(V,c)$ be the vector space spanned by primitive elements in $T\left(V,c\right)$ of homogeneous degree at lest two. 
Let $\mathcal{S}$ be the class of those braided vector spaces $\left(V,c\right)$ whose Nichols algebra $\mathcal{B}\left(V,c\right)$ is obtained dividing out $T\left(V,c\right)$ by its two-sided ideal generated by $E(V,c)$. In Theorem \ref{teo:Ardi}, we also prove that $U(V,c,b)=\mathbb{U}(V,c,\beta)$ whenever $(V,c)$ is in $\mathcal{S}$ (therefore, the main theorem of the present paper extends \cite[Theorem 5.7]{Ardizzoni-Universal}). It is remarkable that, see \cite{Ardizzoni-Sdeg}, this happens, for instance, when
\begin{itemize}
\item $c^{2}=\mathrm{Id}$ and $\mathrm{char}K=0$;
\item there exists $q\in K$, not a root of unity, such that $\left(c+\mathrm{Id}_{V^{\otimes2}}\right)\circ\left(c-q\mathrm{Id}_{V^{\otimes2}}\right)=0$
i.e. $c$ is of Hecke type with mark $q$;
\item the Nichols algebra $\mathcal{B}\left(V,c\right)$ is quadratic.
\item $\left(V,c\right)$ is a braided vector space of diagonal type such
that $\mathcal{B}\left(V,c\right)$ is a domain and its Gelfand-Kirillov
dimension is finite. 
\end{itemize}

Section \ref{sec:mainresult} encloses the main results of the paper. In Proposition \ref{pro: alg is Lie}, we prove that any braided algebra $A$ becomes a braided Lie algebra in a canonical way. In Theorem \ref{thm:Infinitesimal} we show how this Lie algebra structure restricts to the space $P=P(A)$ of primitive elements in $A$ when the braided algebra is further a braided bialgebra. This braided Lie algebra on $P$ is called the infinitesimal braided Lie algebra of $A$ (it is constructed as explained in the first part of the introduction). In Theorem \ref{teo: univ U} a universal property for the universal enveloping algebra of a braided Lie algebra is obtained. Theorems \ref{teo:infpartU} and \ref{teo: generated}, constitutes together a Milnor-Moore type theorem for primitively generated braided bialgebras. Explicitly, Theorem \ref{teo:infpartU} establish that the braided Lie algebra $(V,c,b)$ identifies with the infinitesimal braided Lie algebra of $U(V,c,b)$. Theorem \ref{teo: generated}, which is the main result of the paper, states that every primitively generated braided bialgebra
is isomorphic as a braided bialgebra to the universal enveloping algebra
of its infinitesimal braided Lie algebra. In Example \ref{ex:classic}, we show how to recover the ordinary notion of Lie algebra from that of braided Lie algebra. 

Section \ref{sec:trivbracket}, regards the concept of trivial bracket. In Example \ref{ex: symmetric algebra}, we associates the so-called trivial bracket $b_{tr}$ to any braided vector space $(V,c)$. This way $(V,c,b_{tr})$ becomes a braided Lie algebra and $U(V,c,b_{tr})$ is nothing but the Nichols algebra $\mathcal{B}\left(V,c\right)$ associated to $(V,c)$. In Proposition \ref{pro:vanish}, we establish a condition which guarantees that the bracket of a braided Lie algebra is trivial. In Example \ref{ex:combRank2}, we apply this result to give an example of a braided vector space $\left(V,c\right)$
which becomes a braided Lie algebra only through the trivial bracket. It is remarkable
that this braided vector space does not belong to the class $\mathcal{S}$ mentioned above.

\section{Preliminaries}\label{sec:prelim}

Throughout this paper $K$ will denote a field. All vector spaces
will be defined over $K$ and the tensor product over $K$ will be
denoted by $\otimes$.\medskip{}
 \\
 In this section we define the main notions that we will deal with
in the paper. 
\begin{defn}
Let $V$ be a vector space over a field $K$. A $K$-linear map $c=c_{V}:V\otimes V\rightarrow V\otimes V$
is called a \textbf{braiding} if it satisfies the quantum Yang-Baxter
equation\begin{equation}
c_{1}c_{2}c_{1}=c_{2}c_{1}c_{2}\label{ec: braided equation}\end{equation}
on $V\otimes V\otimes V$, where we set $c_{1}:=c\otimes V$ and $c_{2}:=V\otimes c.$
The pair $\left(V,c\right)$ will be called a\textbf{\ braided vector
space}. A morphism of braided vector spaces $(V,c_{V})$ and $(W,c_{W})$
is a $K$-linear map $f:V\rightarrow W$ such that $c_{W}(f\otimes f)=(f\otimes f)c_{V}.$
\end{defn}
A general method for producing braided vector spaces is to take an
arbitrary braided category $\mathcal{M}$ which is a monoidal subcategory
of the category of $K$-vector spaces. Hence any object $V\in\mathcal{M}$
can be regarded as a braided vector space with respect to $c:=c_{V,V}$,
where $c_{X,Y}:X\otimes Y\rightarrow Y\otimes X$ denotes the braiding
in $\mathcal{M}$, for all $X,Y\in\mathcal{M}$. Both the category
of comodules over a coquasitriangular Hopf algebra and the category
of Yetter-Drinfeld modules over a Hopf algebra with bijective antipode
are examples of such categories. More particularly, every bicharacter
of a group $G$ induces a braiding on the category of $G$-graded
vector spaces.
\begin{defn}
\cite{Ba} A quadruple $(A,m,u,c)$ is called a \textbf{braided algebra}
if
\begin{itemize}
\item $(A,m,u)$ is an associative unital algebra; 
\item $(A,c)$ is a braided vector space;
\item $m$ and $u$ commute with $c$, that is the following conditions
hold: \begin{gather}
c(m\otimes A)=(A\otimes m)(c\otimes A)(A\otimes c),\label{Br2}\\
c(A\otimes m)=(m\otimes A)\left(A\otimes c\right)(c\otimes A),\label{Br3}\\
c(u\otimes A)=A\otimes u,\qquad c(A\otimes u)=u\otimes A.\label{Br4}\end{gather}

\end{itemize}
A morphism of braided algebras is, by definition, a morphism of ordinary
algebras which, in addition, is a morphism of braided vector spaces.

A quadruple $(C,\Delta,\varepsilon,c)$ is called a \textbf{braided
coalgebra} if 
\begin{itemize}
\item $(C,\Delta,\varepsilon)$ is a coassociative counital coalgebra;
\item $(C,c)$ is a braided vector space;
\item $\Delta$ and $\varepsilon$ commute with $c$, that is the following
relations hold: \begin{gather}
(\Delta\otimes C)c=(C\otimes c)(c\otimes C)(C\otimes\Delta),\label{Br5}\\
(C\otimes\Delta)c=(c\otimes C)(C\otimes c)(\Delta\otimes C),\label{Br6}\\
(\varepsilon\otimes C)c=C\otimes\varepsilon,\qquad(C\otimes\varepsilon)c=\varepsilon\otimes C.\label{Br7}\end{gather}

\end{itemize}
A morphism of braided coalgebras is, by definition, a morphism of
ordinary coalgebras which, in addition, is a morphism of braided vector
spaces.

\cite[Definition 5.1]{Ta} A sextuple $(B,m,u,\Delta,\varepsilon,c)$
is a called a \textbf{braided bialgebra} if
\begin{itemize}
\item $(B,m,u,c)$ is a braided algebra
\item $(B,\Delta,\varepsilon,c)$ is a braided coalgebra
\item the following relations hold:\begin{equation}
\Delta m=(m\otimes m)(B\otimes c\otimes B)(\Delta\otimes\Delta).\label{Br1}\end{equation}

\end{itemize}
\end{defn}
Examples of the notions above are algebras, coalgebras and bialgebras
in any braided category $\mathcal{M}$ which is a monoidal subcategory
of the category of $K$-vector spaces. The notion of braided bialgebra
admits a graded counterpart which is called graded braided bialgebra.
For further results on this topic the reader is refereed to \cite{Ardizzoni-Sdeg}.
\begin{example}
\label{ex: tensor algebra} Let $\left(V,c\right)$ be a braided vector
space. Consider the tensor algebra $T=T\left(V\right)$ and let $m_{T}$
and $u_{T}$ denote its multiplication and unit respectively. This
is a graded braided algebra with $n$-th graded component $T^{n}(V)=V^{\otimes n}$.
The braiding $c_{T}$ on $T$ is defined using the the braiding of
$V:$
\end{example}
\begin{figure}\includegraphics[width=3cm,height=2.5cm]{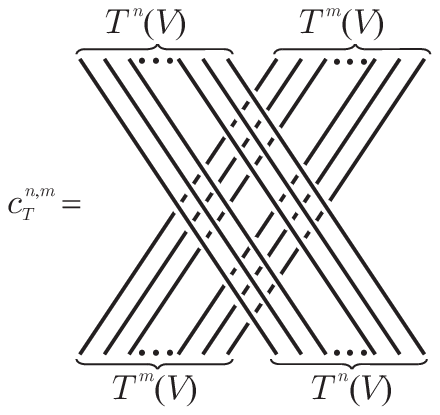} 
\end{figure}

Now $T\otimes T$ becomes itself an algebra with multiplication $m_{T\otimes T}:=\left(m_{T}\otimes m_{T}\right)\circ\left(T\otimes c_{T}\otimes T\right).$
This algebra is denoted by $T\otimes_{c}T.$ The universal property
of the tensor algebra yields two algebra homomorphisms $\Delta_{T}:T\rightarrow T\otimes_{c}T$
and $\varepsilon_{T}:T\rightarrow K$. It is straightforward to check
that $(T,m_{T},u_{T},\Delta_{T},\varepsilon_{T},c_{T})$ is a graded
braided bialgebra.
\begin{defn}
The graded braided bialgebra described in Example \ref{ex: tensor algebra}
is called the \textbf{braided tensor algebra} and is denoted by $T(V,c)$.
\end{defn}
Note that $\Delta_{T}$ really depends on $c$. For instance, one
has $\Delta_{T}\left(z\right)=z\otimes1+1\otimes z+\left(c+\mathrm{Id}\right)\left(z\right)$,
for all $z\in V\otimes V$.
\begin{claim}
\label{nn: connected} Recall that a coalgebra $C$ is called \textbf{connected}
if the coradical $C_{0}$ of $C$ (i.e the sum of all simple subcoalgebras
of $C$) is one dimensional. In this case there is a unique group-like
element $1_{C}\in C$ such that $C_{0}=K1_{C}$. A morphism of connected
coalgebras is just a coalgebra homomorphisms (clearly it preserves
the grouplike element).

By definition, a braided coalgebra $\left(C,c\right)$ is \textbf{connected}
if the underlying coalgebra is connected and, for any $x\in C$, $c(x\otimes1_{C})=1_{C}\otimes x$
and $c(1_{C}\otimes x)=x\otimes1_{C}$

Given a coalgebra $C$ with a distinguished grouplike element $1_{C}$,
we will use the following notation \[
\delta=\delta_{C}:C\rightarrow C\otimes C,\qquad\delta\left(c\right):=c\otimes1_{C}+1_{C}\otimes c-\Delta\left(c\right)\quad\text{for all }c\in C.\]
\end{claim}
\begin{rem}
\label{re: graded connected} Let $C=\bigoplus_{n\in\mathbb{N}}C^{n}$
be a graded braided coalgebra. By \cite[Proposition 11.1.1]{Sw},
if $(C_{n})_{n\in\mathbb{N}}$ is the coradical filtration, then $C_{n}\subseteq\bigoplus_{m\leq n}C^{m}$.
Therefore, $C$ is connected if if $C^{0}$ is a one dimensional vector
space.\end{rem}

\begin{defn}
A braided bialgebra $B$ is called \textbf{primitively
generated} if it is generated as an algebra by the space $P\left(B\right)=\mathrm{Ker}\delta_{B}=\left\{ b\in B\mid\Delta\left(b\right)=b\otimes1_{B}+1_{B}\otimes b\right\} $
of primitive elements in $B$. See \cite[page 239]{Milnor-Moore}.\end{defn}

\begin{example}\label{ex: cocommutative}
 Let $B$ be a connected ordinary bialgebra in characteristic zero. 

Assume that $B$ is primitively generated and set $P=P(A)$. Then, we have a Hopf algebra projection $T(P)\rightarrow B$. Since $T(P)$ is cocommutative, $B$ is cocommutative too. 

Conversely, if $B$ is cocommutative, then $\mathrm{gr}B$ is also cocommutative connected. By \cite[$(3)\Rightarrow(4)$ in Theorem 2.15]{AMS-MM2}, we have that $\mathrm{gr}B$ is strictly generated as an algebra by $\mathrm{gr}^1B=B_1/B_0$. Hence $B$ is primitively generated. Note that we did not used Milnor-Moore theorem to conclude that $B$ is primitively generated. 
\end{example}

\begin{prop}
\label{pro:primgen-is-connected}\cite[Proposition 5.8]{Ardizzoni-Universal}
Let $B$ be a primitively generated braided bialgebra.
Then the underlying braided coalgebra is connected.\end{prop}
\begin{example}
We give an example of a finite dimensional connected braided bialgebra
which is not primitively generated (it is a particular case of the
graded dual of $S$ as in \cite[Example 2.5]{AS}). Assume $\mathrm{char}K=2$. Consider the algebra $B$ generated
as an algebra over $K$ by $x$ and $y$ with relations $x^{2},y^{2},yx-xy$.
Then $1,x,y,xy$ form a basis for $B$. We can regard $B$ as an ordinary
bialgebra by setting \[
\Delta_{B}\left(x\right)=1\otimes x+x\otimes1,\qquad\Delta_{B}\left(y\right)=1\otimes y+x\otimes x+y\otimes1.\]
Clearly $B$ is not primitively generated. On the other hand, if we
set $B\left(0\right):=K,B\left(1\right)=K+Kx,B\left(2\right):=K+Kx+Ky,$
and $B\left(3\right)=B$, we get a filtration for the coalgebra $B$.
By \cite[Proposition 11.1.1]{Sw}, this yields $\mathrm{Corad}\left(B\right)\subseteq B\left(0\right)=K$.
Hence $B$ is connected.\end{example}
\begin{rem}
It is remarkable that, in characteristic zero, in the literature
there is no example of a finite dimensional connected braided bialgebra
which is not primitively generated. See also \cite[Conjecture 1.4]{AS- FiniteQuantCartan}
and \cite[Conjecture 2.7]{AS}.\end{rem}
\begin{lem}
\label{lem: lim of br bialg}Let $\left(\left(B_{u}\right)_{u\in\mathbb{N}},\left(\xi_{u}^{v}\right)_{u,v\in\mathbb{N}}\right)$
be a direct system of vector spaces, where, for $u\leq v,$ $\xi_{u}^{v}:B_{u}\rightarrow B_{v}.$
Assume that each $B_{u}$ is endowed with a braided bialgebra structure
such that $\xi_{u}^{v}$ is a braided bialgebra homomorphism, for
every $u,v\in\mathbb{N}.$ Then $\underrightarrow{\lim}B_{u}$ carries
a natural braided bialgebra structure that makes it the direct limit
of $\left(\left(B_{u}\right)_{u\in\mathbb{N}},\left(\xi_{u}^{v}\right)_{u,v\in\mathbb{N}}\right)$
as a direct system of braided bialgebras. \end{lem}
\begin{proof}
It is straightforward. %

\end{proof}

\section{Pre-categorical and categorical subspaces}\label{sec:pre-categ}
In this section we investigate the properties of pre-categorical and categorical vector spaces that will be needed next section.

\begin{defn}
\label{def: categorical}(cf. \cite[2.2]{Kharchenko- connected}) A subspace
$L$ of a braided vector space $\left(V,c\right)$ will be called
\textbf{pre-categorical} whenever\begin{equation}
c\left(L\otimes V+V\otimes L\right)\subseteq L\otimes V+V\otimes L.\label{eq:precategorical}\end{equation}
\end{defn}

\begin{prop}\label{pro:quotPreCateg}
Let $\left(V,c\right)$ be a braided vector space.

1) Let $f:\left(V,c\right)\rightarrow\left(V',c'\right)$ be a morphism
of braided vector spaces. Then, $\mathrm{Ker}\left(f\right)$ is a
pre-categorical subspace of $(V,c)$.

2) Let $L$ be pre-categorical subspace of $\left(V,c\right)$.
Then, the quotient $V/L$ carries a unique braided vector
space structure such that the canonical projection is a morphism of
braided vector spaces. \end{prop}
\begin{proof}
1) We have\[
\left(f\otimes f\right)c\left(\mathrm{Ker}\left(f\otimes f\right)\right)=c'\left(f\otimes f\right)\left(\mathrm{Ker}\left(f\otimes f\right)\right)=0\]
so that $c\left(\mathrm{Ker}\left(f\otimes f\right)\right)\subseteq\mathrm{Ker}\left(f\otimes f\right)$.
Since $\mathrm{Ker}\left(f\otimes f\right)=\mathrm{Ker}\left(f\right)\otimes V+V\otimes\mathrm{Ker}\left(f\right)$,
we get that $\mathrm{Ker}\left(f\right)$ is a pre-categorical subspace
of $(V,c)$.

2) Set $V':=V/L$ and denote by $\pi:V\rightarrow V'$ the canonical
projection. We have\[
c\left(\mathrm{Ker}\left(\pi\otimes\pi\right)\right)=c\left(L\otimes V+V\otimes L\right)\overset{\eqref{eq:precategorical}}{\subseteq}L\otimes V+V\otimes L=\mathrm{Ker}\left(\pi\otimes\pi\right)\]
so that there exists $c':V'\otimes V'\rightarrow V'\otimes V'$ such
that $c'\left(\pi\otimes\pi\right)=\left(\pi\otimes\pi\right)c$.
By this equality and the fact that $\pi$ is subjective, it is straightforward
to prove that $\left(V',c'\right)$ is a braided vector space.\end{proof}
\begin{defn}
A subspace $L$ of a braided vector space $\left(V,c\right)$ is called
\textbf{categorical} (cf. \cite[Section 6]{Ta}; see also \cite[2.5]{AndrusNatale-BraidedHopfMatched})
whenever \begin{equation}
c\left(L\otimes V\right)\subseteq V\otimes L\qquad\text{and}\qquad c\left(V\otimes L\right)\subseteq L\otimes V.\label{form:categorical}\end{equation}
 In this case $c$ induces maps $c_{L,V}:L\otimes V\rightarrow V\otimes L$
and $c_{V,L}:V\otimes L\rightarrow L\otimes V.$

If $L$ and $L^{\prime}$ are categorical subspaces of $\left(V,c\right)$,
we say that $f:L\rightarrow L'$ is a \textbf{morphism of categorical
subspaces} of $\left(V,c\right)$ if $f$ is a $K$-linear map such
that \begin{equation}
c_{L',V}\left(f\otimes V\right)=\left(V\otimes f\right)c_{L,V}\qquad\text{and}\qquad c_{V,L'}\left(V\otimes f\right)=\left(f\otimes V\right)c_{V,L}.\label{form:betaextended}\end{equation}

Let $L$ and $L^{\prime}$ be categorical subspaces of $\left(V,c\right).$
Then\[
c\left(L\otimes L^{\prime}\right)\subseteq c\left(L\otimes V\right)\cap c\left(V\otimes L^{\prime}\right)\subseteq\left(V\otimes L\right)\cap\left(L^{\prime}\otimes V\right)=L^{\prime}\otimes L\]
 so that $c$ induces a map $c_{L,L^{\prime}}:L\otimes L^{\prime}\rightarrow L^{\prime}\otimes L.$
If in addition $L\subseteq L'$, we get that $L$ is a categorical
subspace of $L'$ too. 

We will sometime use the simplified notation $c_{L}:=c_{L,L}.$ \end{defn}
\begin{rem}
\label{rem:precategorical}Any categorical subspace of a braided vector
space is pre-categorical.
\end{rem}
In the connected case, the following result appeared in \cite[page 4]{Kharchenko- connected}.
\begin{lem}
\label{lem: P categorical}\cite[Lemma 2.10]{Ardizzoni-Universal}
Let $\left(B,c\right)$ be a braided bialgebra. Let $P=P\left(B\right)$
be the space of primitive elements of $B.$ Then $P$ is a categorical
subspace of $\left(B,c\right)$. \end{lem}
\begin{defn}
\label{infinitesimal} With same assumptions and notations as in Lemma
\ref{lem: P categorical}, the braiding $c_{P}$ of $P$ (which is
induced by $c$) will be called the \textbf{infinitesimal braiding
of} $B$ while $(P,c_{P})$ will be called the \textbf{infinitesimal
part of} $B$.\end{defn}
\begin{lem}
\label{lem:morphcat}Let $L$ and $L^{\prime}$ be categorical subspaces
of $\left(V,c\right)$, and let $f:L\rightarrow L'$ be a morphism
of categorical subspaces of $\left(V,c\right)$. Then $\mathrm{Ker}\left(f\right)$
and $\mathrm{Im}\left(f\right)$ are categorical subspaces of $\left(V,c\right)$.\end{lem}
\begin{proof}
Set $K:=\mathrm{Ker}\left(f\right)$. Then\[
\left(V\otimes f\right)c\left(K\otimes V\right)=\left(V\otimes f\right)c_{L,V}\left(K\otimes V\right)\overset{\eqref{form:betaextended}}{=}c_{L',V}\left(f\otimes V\right)\left(K\otimes V\right)=0\]
so that $c\left(K\otimes V\right)\subseteq\mathrm{Ker}\left(V\otimes f\right)=V\otimes K$.
Similarly one gets $c\left(V\otimes K\right)\subseteq K\otimes V$
so that $K$ is a categorical subspaces of $\left(V,c\right)$. Set
$D:=\mathrm{Im}\left(f\right)$. We have\[
c\left(D\otimes V\right)=c_{L',V}\left(D\otimes V\right)=c_{L',V}\left(f\otimes V\right)\left(L\otimes V\right)\overset{\eqref{form:betaextended}}{=}\left(V\otimes f\right)c_{L,V}\left(L\otimes V\right)\subseteq V\otimes D\]
so that $c\left(D\otimes V\right)\subseteq V\otimes D$. Similarly
one gets $c\left(V\otimes D\right)\subseteq D\otimes V$ so that $D$
is a categorical subspaces of $\left(V,c\right)$. \end{proof}
\begin{defn}
\label{def: powers of V}Let $B$ be a braided bialgebra. Let $V$
be a vector subspace of $B.$ Then we can define recursively the iterated
powers of $V$ in $B$ by setting \[
V^{\cdot_{B}0}:=K\qquad\text{and}\qquad V^{\cdot_{B}n}:=V^{\cdot_{B}\left(n-1\right)}\cdot_{B}V,\qquad\text{for all }n\geq1.\]
\end{defn}
\begin{lem}
\label{lem: powers categorical}Let $\left(B,c\right)$ be a braided
bialgebra. Let $V$ be a categorical subspace of $\left(B,c\right).$
Then $V^{\cdot_{B}n}$ is a categorical subspace of $B$ for all $n\in\mathbb{N}.$\end{lem}
\begin{proof}
We prove the statement by induction on $n\in\mathbb{N}$. Set $V^{n}:=V^{\cdot_{B}n}$.
For $n=0$ it follows by (\ref{Br4}) and for $n=1$ there is nothing
to prove. Let $n\geq2$ and assume $V^{\cdot_{B}n-1}$ is a categorical
subspace of $\left(B,c\right)$. Let $m$ denote the multiplication
of $B$. Then\begin{eqnarray*}
c\left(V^{n}\otimes B\right) & \subseteq & c\left(V^{n-1}\cdot_{B}V\otimes B\right)\\
 & \subseteq & c\left(m\otimes B\right)\left(V^{n-1}\otimes V\otimes B\right)\\
 & \overset{\left(\ref{Br2}\right)}{\subseteq} & (B\otimes m)(c\otimes B)(B\otimes c)\left(V^{n-1}\otimes V\otimes B\right)\\
 & \subseteq & (B\otimes m)(c\otimes B)\left(V^{n-1}\otimes B\otimes V\right)\\
 & \overset{\text{induction hypothesis}}{\subseteq} & (B\otimes m)\left(B\otimes V^{n-1}\otimes V\right)\subseteq B\otimes V^{n}.\end{eqnarray*}
 Similarly, using (\ref{Br3}) one proves that $c\left(B\otimes V^{n}\right)\subseteq V^{n}\otimes B.$ \end{proof}
\begin{prop}
\label{pro: cat-braid}Let $B$ be a braided bialgebra
generated as an algebra by some subspace $V$. Then, the following
assertions are equivalent. 
\begin{enumerate}
\item [$\left(i\right)$] $V$ is a braided subspace of $B.$ 
\item [$\left(ii\right)$] $V$ is a categorical subspace of $B.$ 
\end{enumerate}
\end{prop}
\begin{proof} Denote by $c$ the braiding of $B$.

$\left(ii\right)\Rightarrow\left(i\right)$ By Definition \ref{def: categorical},
$c$ induces a map $c_{V}:V\otimes V\rightarrow V\otimes V.$ Clearly
$\left(V,c_{V}\right)$ is a braided subspace of $\left(B,c\right).$

$\left(i\right)\Rightarrow\left(ii\right)$ We have to check that
$c\left(V\otimes B\right)\subseteq B\otimes V$ and $c\left(B\otimes V\right)\subseteq V\otimes B.$
We treat the former inclusion, the latter having a similar proof.
Since $B$ is generated by $V$, then $B=\cup V^{n}$ where $V^{n}:=V^{\cdot_{B}n}$
for all $n\in\mathbb{N}.$ Thus it suffices to check that $c\left(V\otimes V^{n}\right)\subseteq V^{n}\otimes V.$
For $n=1,$ there is nothing to prove as $V$ is a braided subspace
of $\left(B,c\right)$. Assume that the formulas are true for $n-1$
and denote by $m$ the multiplication of $B$. Then, we have \begin{eqnarray*}
c\left(V\otimes V^{n}\right) & = & c\left(B\otimes m\right)\left(V\otimes V^{n-1}\otimes V\right)\\
 & \overset{(\ref{Br3})}{=} & (m\otimes B)\left(B\otimes c\right)(c\otimes B)\left(V\otimes V^{n-1}\otimes V\right)\\
 & \overset{\text{induction hypothesis}}{\subseteq} & (m\otimes B)\left(B\otimes c\right)\left(V^{n-1}\otimes V\otimes V\right)\\
 & \subseteq & (m\otimes B)\left(V^{n-1}\otimes V\otimes V\right)=V^{n}\otimes V.\end{eqnarray*}

\end{proof}

\section{Braided Ideals}\label{sec:brideal}

In this section we face the problem of quotienting braided bialgebras. In particular, we investigate braided bialgebra quotients of the braided tensor algebra $T(V,c)$ associated to a given braided vector space $(V,c)$. These braided bialgebras will be characterized as $(V,c)$-bialgebras.

\begin{defn}
A \textbf{braided ideal} of a braided bialgebra $B$
is a pre-categorical subspace of $B$ which is
both an ideal and a coideal of $B$.\end{defn}

\begin{lem}
\label{pro: braided quotient}For a given braided braided ideal $I$
of a braided bialgebra $B,$ the quotient $B/I$ carries a unique
braided bialgebra structure such that the canonical projection is
a braided bialgebra homomorphism. \end{lem}
\begin{proof}
Set $A:=B/I$. Since $I$ is both an ideal and a coideal of $B$,
it is clear that $A$ is both an algebra and a coalgebra in such a
way that the projection $\pi:B\rightarrow A$ is both an algebra and
a coalgebra map.  Since $I$ is be pre-categorical subspace of $B$, by Proposition \ref{pro:quotPreCateg}, the quotient $A=B/I$ carries a unique braided vector
space structure such that $\pi$ is a morphism of
braided vector spaces. From this fact, since $A$ is a braided bialgebra, $\pi$ is surjective
and it is both an algebra and a coalgebra map, it is straightforward
to prove that $B$ is a braided bialgebra in such a way that $\pi$
is a braided bialgebra map. 
\end{proof}
The following result in the case $B=T\left(V,c\right)$ for some braided
vector space $\left(V,c\right)$ appeared in \cite[page 7]{Kharchenko- connected}.
\begin{thm}
\label{thm:idealprecategorical}Let $B$ be a braided bialgebra with
infinitesimal part $\left(P,c_{P}\right)$. Let $S$ be a pre-categorical
subspace of $\left(P,c_{P}\right)$. Then the ideal $\left(S\right)$
of $B$ generated by $S$ is a braided ideal of $B$.\end{thm}
\begin{proof}
Let $c,m,\Delta$ and $\varepsilon$ denote the braiding, the multiplication
the comultiplication and the counit of $B$ respectively. Set $I:=\left(S\right)=BSB.$
We prove that $c\left(I\otimes B\right)\in I\otimes B+B\otimes I.$
We have\begin{eqnarray*}
c\left(I\otimes B\right) & = & c\left(BSB\otimes B\right)\\
 & = & c\left(m\otimes B\right)\left(B\otimes m\otimes B\right)\left(B\otimes S\otimes B\otimes B\right)\\
 & = & \left(B\otimes m\right)\left(B\otimes B\otimes m\right)\left(c\otimes B\otimes B\right)\left(B\otimes c\otimes B\right)\left(B\otimes S\otimes c\left(B\otimes B\right)\right)\\
 & \subseteq & \left(B\otimes m\right)\left(B\otimes B\otimes m\right)\left(c\otimes B\otimes B\right)\left(B\otimes c\otimes B\right)\left(B\otimes S\otimes B\otimes B\right)\\
 & \overset{\text{\eqref{eq:precategorical}}}{\subseteq} & \left(B\otimes m\right)\left(B\otimes B\otimes m\right)\left(c\otimes B\otimes B\right)\left(B\otimes S\otimes B\otimes B+B\otimes B\otimes S\otimes B\right)\\
 & \overset{\text{\eqref{eq:precategorical}}}{\subseteq} & \left(B\otimes m\right)\left(B\otimes B\otimes m\right)\left(B\otimes S\otimes B\otimes B+S\otimes B\otimes B\otimes B+B\otimes B\otimes S\otimes B\right)\\
 & \subseteq & B\otimes SBB+S\otimes BBB+B\otimes BSB\subseteq B\otimes I+I\otimes B.\end{eqnarray*}
Similarly one gets $c\left(B\otimes I\right)\subseteq I\otimes B+B\otimes I$
so that $I$ is a pre-categorical subspace of $\left(B,c\right)$.
It remains to check that $I$ is a coideal of $B$. Note that $S\subseteq P$
so that $\Delta\left(S\right)\subseteq B\otimes S+S\otimes B$ and
$\varepsilon\left(S\right)=0.$ Let us check that $\Delta\left(I\right)\in I\otimes B+B\otimes I:$\begin{eqnarray*}
 &  & \Delta\left(I\right)\\
 & = & \Delta\left(SBS\right)\\
 & = & \Delta m\left(m\otimes B\right)\left(B\otimes S\otimes B\right)\\
 & \overset{\text{(\ref{Br1})}}{=} & \left(m\otimes m\right)\left(B\otimes c\otimes B\right)\left(\Delta\otimes\Delta\right)\left(m\otimes B\right)\left(B\otimes S\otimes B\right)\\
 & \overset{\text{(\ref{Br1})}}{=} & \left[\begin{array}{c}
\left(m\otimes m\right)\left(B\otimes c\otimes B\right)\left(m\otimes m\otimes B\otimes B\right)\\
\left(B\otimes c\otimes B\otimes B\otimes B\right)\left(\Delta\otimes\Delta\otimes\Delta\right)\left(B\otimes S\otimes B\right)\end{array}\right]\\
 & \subseteq & \left[\begin{array}{c}
\left(m\otimes m\right)\left(B\otimes c\otimes B\right)\left(m\otimes m\otimes B\otimes B\right)\\
\left(B\otimes c\otimes B\otimes B\otimes B\right)\left(B\otimes B\otimes B\otimes S\otimes B\otimes B+B\otimes B\otimes S\otimes B\otimes B\otimes B\right)\end{array}\right]\\
 & \overset{\text{(\ref{eq:precategorical})}}{\subseteq} & \left[\begin{array}{c}
\left(m\otimes m\right)\left(B\otimes c\otimes B\right)\left(m\otimes m\otimes B\otimes B\right)\\
\left(B\otimes B\otimes B\otimes S\otimes B\otimes B+B\otimes B\otimes S\otimes B\otimes B\otimes B+B\otimes S\otimes B\otimes B\otimes B\otimes B\right)\end{array}\right]\\
 & \subseteq & \left(m\otimes m\right)\left(B\otimes c\otimes B\right)\left(B\otimes BS\otimes B\otimes B+B\otimes SB\otimes B\otimes B+BS\otimes B\otimes B\otimes B\right)\\
 & \subseteq & \left(m\otimes m\right)\left(B\otimes c\otimes B\right)\left(B\otimes I\otimes B\otimes B+I\otimes B\otimes B\otimes B\right)\\
 & \subseteq & \left(m\otimes m\right)\left(B\otimes I\otimes B\otimes B+B\otimes B\otimes I\otimes B+I\otimes B\otimes B\otimes B\right)\\
 & \subseteq & BI\otimes B+B\otimes IB+IB\otimes B\\
 & \subseteq & I\otimes B+B\otimes I.\end{eqnarray*}
 Furthermore $\varepsilon\left(I\right)=\varepsilon\left(BSB\right)=\varepsilon\left(B\right)\varepsilon\left(S\right)\varepsilon\left(B\right)=0.$
Therefore $I$ is a coideal.\end{proof}
\begin{cor}
\label{coro:idealcategorical}Let $B$ be a braided bialgebra with
infinitesimal part $\left(P,c_{P}\right)$. Let $S$ be a categorical
subspace of $\left(P,c_{P}\right)$. Then the ideal $\left(S\right)$
of $B$ generated by $S$ is a braided ideal of $B$.\end{cor}
\begin{proof}
Apply Remark \ref{rem:precategorical} and Theorem \ref{thm:idealprecategorical}. \end{proof}
\begin{prop}
\label{pro: W}Let $B$ be a braided bialgebra with infinitesimal
part $P$. Let $b:P\rightarrow P$ be a morphism of categorical subspaces
of $B$. Then the two-sided ideal of $B$ generated by $\left[\mathrm{Id}-b\right]\left[P\right]$
is a braided ideal. Thus the quotient \[
W\left(B,b\right):=\frac{B}{\left(\left[\mathrm{Id}-b\right]\left[P\right]\right)}\]
carries a unique braided bialgebra structure such that the canonical
projection is a braided bialgebra homomorphism. \end{prop}
\begin{proof}
Since both $\mathrm{Id}_{P}:P\rightarrow P$ and $b:P\rightarrow P$
are morphisms of categorical subspaces of $B$, so is their difference
$\mathrm{Id}-b:P\rightarrow P$. By Lemma \ref{lem:morphcat}, $S:=\left[\mathrm{Id}-b\right]\left[P\right]$
is a categorical subspace of $B$ whence of $P$. By Theorem \ref{coro:idealcategorical},
$I=\left(S\right)=\left(\left[\mathrm{Id}-b\right]\left[P\right]\right)$
is a braided ideal of $B$. The conclusion follows by Lemma \ref{pro: braided quotient}.\end{proof}
\begin{lem}
\label{lem: beta extended}Let $B$ be a braided bialgebra with infinitesimal
part $P$. Let $V$ be a categorical subspace of $B$ generating $B$
as an algebra and contained in $P$. Let $j:V\rightarrow P$ be
the canonical inclusion. The following are equivalent for a $K$-linear
map $b:P\rightarrow V:$ 

(i) $jb:P\rightarrow P$ is a morphism of categorical subspaces of
$B$;

(ii) $b:P\rightarrow V$ fulfills\begin{equation}
c_{V}\left(b\otimes V\right)=\left(V\otimes b\right)c_{P,V}\qquad\text{and}\qquad c_{V}\left(V\otimes b\right)=\left(b\otimes V\right)c_{V,P}.\label{form: beta}\end{equation}
\end{lem}
\begin{proof}
$(i)\Rightarrow\left(ii\right)$ Clearly $b$ fulfills (\ref{form: beta})
if $jb$ fulfills (\ref{form:betaextended}). 

$(ii)\Rightarrow\left(i\right)$ Since $B$ is generated by $V$,
we have that $B=\cup V^{n}$ where $V^{n}:=V^{\cdot_{B}n}$ for all
$n\in\mathbb{N}.$ Hence, in order to prove that $jb$ fulfills (\ref{form:betaextended})
it suffices to check, by induction on $n\geq1,$ that $b$ fulfills
\[
c_{V,V^{n}}\left(b\otimes V^{n}\right)=\left(V^{n}\otimes b\right)c_{P,V^{n}}\qquad c_{V^{n},V}\left(V^{n}\otimes b\right)=\left(b\otimes V^{n}\right)c_{V^{n},P}.\]
 Note that the notations above make sense in view of Lemma \ref{lem: P categorical}
and Lemma \ref{lem: powers categorical}. For $n=1$ there is nothing
to prove as $c_{V,V}=c_{V}$. Assume that the formulas are true for
$n-1.$ Then, we have \begin{eqnarray*}
 &  & c_{V,V^{n}}\left(b\otimes V^{n}\right)\left(P\otimes m_{V^{n-1},V}\right)\\
 & = & c_{V,V^{n}}\left(V\otimes m_{V^{n-1},V}\right)\left(b\otimes V^{n-1}\otimes V\right)\\
 & \overset{(\ref{Br2})}{=} & \left(m_{V^{n-1},V}\otimes V\right)\left(V^{n-1}\otimes c_{V}\right)\left(c_{V,V^{n-1}}\otimes V\right)\left(b\otimes V^{n-1}\otimes V\right)\\
 & = & \left(m_{V^{n-1},V}\otimes V\right)\left(V^{n-1}\otimes c_{V}\right)\left(V^{n-1}\otimes b\otimes V\right)\left(c_{P,V^{n-1}}\otimes V\right)\\
 & = & \left(m_{V^{n-1},V}\otimes V\right)\left(V^{n-1}\otimes V\otimes b\right)\left(V^{n-1}\otimes c_{P,V}\right)\left(c_{P,V^{n-1}}\otimes V\right)\\
 & = & \left(V^{n}\otimes b\right)\left(m_{V^{n-1},V}\otimes P\right)\left(V^{n-1}\otimes c_{P,V}\right)\left(c_{P,V^{n-1}}\otimes V\right)\\
 & \overset{(\ref{Br2})}{=} & \left(V^{n}\otimes b\right)c_{P,V^{n}}\left(P\otimes m_{V^{n-1},V}\right)\end{eqnarray*}
 where $m_{V^{n-1},V}:V^{n-1}\otimes V\rightarrow V^{n}$ is induced
by the multiplication of $B$. Since $P\otimes m_{V^{n-1},V}$ is
surjective, we arrive at $c_{V,V^{n}}\left(b\otimes V^{n}\right)=\left(V^{n}\otimes b\right)c_{P,V^{n}}.$
Similarly one gets $c_{V^{n},V}\left(V^{n}\otimes b\right)=\left(b\otimes V^{n}\right)c_{V^{n},P}.$ \end{proof}
\begin{cor}
\label{coro: W}Let $B$ be a braided bialgebra with infinitesimal
part $P$. Let $V$ be a braided subspace of $B$ generating $B$
as an algebra and contained in $P$. Let $b:P\rightarrow V$ fulfill
(\ref{form: beta}). Then \[
W\left(B,b\right):=\frac{B}{\left(\left[\mathrm{Id}-b\right]\left[P\right]\right)}\]
 carries a unique connected braided bialgebra structure such that
the canonical projection $\pi_{W}:B\rightarrow W\left(B,b\right)$
is a braided bialgebra homomorphism. \end{cor}
\begin{proof}
It follows by Proposition \ref{pro: cat-braid}, Lemma \ref{lem: beta extended}
and Proposition \ref{pro: W}. \end{proof}
\begin{defn}
Let $\left(V,c\right)$ be a braided vector space. A $\left(V,c\right)$\textbf{-bialgebra}
is a pair $\left(B,\gamma\right)$ such that 
\begin{itemize}
\item $B$ is a braided bialgebra with infinitesimal part $\left(P,c_{P}\right)$; 
\item $\gamma:\left(V,c\right)\rightarrow\left(P,c_{P}\right)$ is a morphism
of braided vector spaces; 
\item $\gamma\left(V\right)$ generates $B$ as a $K$-algebra. 
\end{itemize}
A \textbf{morphism} $\pi:\left(B,\gamma\right)\rightarrow\left(B^{\prime},\gamma^{\prime}\right)$
\textbf{of} $\left(V,c\right)$\textbf{-bialgebras} is a braided bialgebra
homomorphism $\pi:B\rightarrow B^{\prime}$ such that $P\left(\pi\right)\gamma=\gamma^{\prime}$
where $P\left(\pi\right):P\left(B\right)\rightarrow P\left(B^{\prime}\right)$
is the obvious induced morphism. 
\end{defn}
\begin{rem}
Let $\left(V,c\right)$ be a braided vector space. Then any $\left(V,c\right)$-bialgebra
$\left(B,\gamma\right)$ is connected. In fact, $\gamma\left(V\right)\subseteq P(B)$
generates $B$ as a $K$-algebra. Hence $B$ is primitively generated.
Thus Proposition \ref{pro:primgen-is-connected} applies.\end{rem}
\begin{lem}\label{lem:quotOfTensor}
Let $\left(V,c\right)$ be a braided vector space. The following assertions
are equivalent for a braided bialgebra $B$.
\begin{enumerate}
\item There is a map $\gamma$ such that $\left(B,\gamma\right)$ is a $\left(V,c\right)$-bialgebra.
\item There is a surjective braided bialgebra map $\varphi:T\left(V,c\right)\rightarrow B$.
\end{enumerate}
\end{lem}
\begin{proof}
Set $T:=T\left(V,c\right)$.

$\left(1\right)\Rightarrow\left(2\right)$ By the universal property
of the braided tensor algebra, there is a unique braided bialgebra
map $\varphi:T\rightarrow B$. Since $\gamma\left(V\right)$ generates
$B$ as a $K$-algebra, then $\varphi$ is surjective.

$\left(2\right)\Rightarrow\left(1\right)$ Since $\varphi$ is a braided
bialgebra map, it preserves primitive elements, so that one has $\varphi\left(V\right)\subseteq\varphi\left(P\left(T\right)\right)\subseteq P\left(\varphi\left(T\right)\right)\subseteq P\left(B\right)$.
Thus we can set $\gamma$ to be the restriction of $\varphi$ to $V$.
Since $\varphi$ is compatible with braidings, one has that $\gamma$
is a morphism of braided vector spaces. Since $\varphi$ is surjective,
then $\gamma\left(V\right)$ generates $B$ as a $K$-algebra. In
other words $\left(B,\gamma\right)$ is a $\left(V,c\right)$-bialgebra.\end{proof}
\begin{thm}
\label{teo: (V,c)-bialg}Let $\left(V,c\right)$ be a braided vector
space and let $\left(B,\gamma\right)$ be a $\left(V,c\right)$-bialgebra
endowed with a map $b:P\left(B\right)\rightarrow V$ such that \begin{eqnarray}
c_{\gamma\left(V\right)}\left(\gamma b\otimes\gamma\left(V\right)\right) & = & \left(\gamma\left(V\right)\otimes\gamma b\right)c_{P\left(B\right),\gamma\left(V\right)},\label{form: bracketable 1}\\
c_{\gamma\left(V\right)}\left(\gamma\left(V\right)\otimes\gamma b\right) & = & \left(\gamma b\otimes\gamma\left(V\right)\right)c_{\gamma\left(V\right),P\left(B\right)}.\label{form: bracketable 2}\end{eqnarray} Then \[
W=W\left(B,\gamma b\right)=\frac{B}{\left(\left[\mathrm{Id}-\gamma b\right]\left[P\left(B\right)\right]\right)}\]
carries a unique $\left(V,c\right)$-bialgebra structure such that
the canonical projection $\pi_{W}:B\rightarrow W$ is a $\left(V,c\right)$-bialgebra
homomorphism. \end{thm}
\begin{proof}
$\gamma\left(V\right)$ is a braided subspace of $B$ generating $B$
as a $K$-algebra and contained in $P\left(B\right).$ In view of
(\ref{form: bracketable 1}) and (\ref{form: bracketable 2}), we
have that $\gamma b$ fulfills (\ref{form: beta}) where $\gamma\left(V\right)$
plays the role of $V.$ By Corollary \ref{coro: W}, $W\left(B,\gamma b\right)$
carries a unique braided bialgebra structure such that the canonical
projection $\pi_{W}:B\rightarrow W$ is a braided bialgebra homomorphism.
Clearly $\pi_{W}$ induces a map $P\left(\pi_{W}\right):P\left(B\right)\rightarrow P\left(W\right).$
Then $\left(W,P\left(\pi_{W}\right)\circ\gamma\right)$ is a $\left(V,c\right)$-bialgebra
and, by construction, $\pi_{W}:B\rightarrow W$ is a $\left(V,c\right)$-bialgebra
homomorphism. \end{proof}
\begin{defn}
\label{def: quotient (V,c)-bialg}With hypothesis and notations of
Theorem \ref{teo: (V,c)-bialg}, $\left(W,P\left(\pi_{W}\right)\circ\gamma\right)$
will be called the \textbf{quotient} $\left(V,c\right)$\textbf{-bialgebra
of} $\left(B,\gamma\right)$ \textbf{induced by} $b.$ 
\end{defn}

\section{Universal enveloping algebra}\label{sec:univenv}

In this section we introduce and investigate the notion of braided Lie algebra for a braided vector space. We start by introducing the concept of bracket of finite rank and the corresponding universal enveloping algebra.

\begin{defn}
\label{def: main construction}Let $\left(V,c\right)$ be a braided
vector space and let $n\in\mathbb{N}$. 

Suppose there are $\left(V,c\right)$-bialgebras
\[
\left(U^{\left[0\right]},i^{\left[0\right]}\right),\left(U^{\left[1\right]},i^{\left[1\right]}\right),\dots,\left(U^{\left[n\right]},i^{\left[n\right]}\right)\]
such that 
\begin{itemize}
\item $U^{\left[0\right]}:=T\left(V,c\right)$ and $i^{\left[0\right]}:V\rightarrow P\left(U^{\left[0\right]}\right)$
is the restriction of the inclusion $V\rightarrow T\left(V,c\right)$;
\item there exists a map $b^{\left[i\right]}:P\left(U^{\left[i\right]}\right)\rightarrow V$ such that \eqref{form: bracketable 1} and \eqref{form: bracketable 2} hold for $(B,\gamma,b)=\left(U^{\left[i\right]},i^{\left[i\right]},b^{\left[i\right]}\right)$, for every $0\leq i\leq n$;
\item $\left(U^{\left[i\right]},i^{\left[i\right]}\right)$ is the quotient
$\left(V,c\right)$-bialgebra of $\left(U^{\left[i-1\right]},i^{\left[i-1\right]}\right)$
induced by $b^{\left[i-1\right]}$ for every $1\leq i\leq n$.
\end{itemize}
In this case we will say that\textbf{ \[
\left(b^{\left[i\right]}\right)_{0\leq i\leq n}\]
}is a \textbf{bracket of rank $n$ for a braided vector space $\left(V,c\right)$}.

Given a bracket $\left(b^{\left[i\right]}\right)_{0\leq i\leq n}$
of rank $n$ for $\left(V,c\right)$ we can define \textbf{the universal
enveloping algebra of rank} $n+1$ as the quotient $\left(V,c\right)$-bialgebra
$\left(U^{\left[n+1\right]},i^{\left[n+1\right]}\right)$ of $\left(U^{\left[n\right]},i^{\left[n\right]}\right)$
induced by $b^{\left[n\right]}$ (see Definition \ref{def: quotient (V,c)-bialg}).
This will be denoted also by \[
U\left(V,c,\left(b^{\left[i\right]}\right)_{0\leq i\leq n}\right).\]
Note that, with this notation, we have that $U\left(V,c,\left(b^{\left[i\right]}\right)_{0\leq i\leq t-1}\right)=\left(U^{\left[t\right]},i^{\left[t\right]}\right)$
for all $1\leq t\leq n$. This leads to define the \textbf{the universal
enveloping algebra of rank} $0$ just as $U\left(V,c\right):=\left(U^{\left[0\right]},i^{\left[0\right]}\right)$.

For sake of simplicity we will adopt the following notation for all
$n\in\mathbb{N}$ 

\[
P^{\left[n\right]}:=P\left(U^{\left[n\right]}\right),\qquad V^{\left[n\right]}:=i^{\left[n\right]}\left(V\right).\]
 Hence for each $n\in\mathbb{N}$, (\ref{form: bracketable 1}) and
(\ref{form: bracketable 2}) take form \begin{eqnarray}
c_{V^{\left[n\right]}}\left(i^{\left[n\right]}b^{\left[n\right]}\otimes V^{\left[n\right]}\right) & = & \left(V^{\left[n\right]}\otimes i^{\left[n\right]}b^{\left[n\right]}\right)c_{P^{\left[n\right]},V^{\left[n\right]}},\label{form:bracket1}\\
c_{V^{\left[n\right]}}\left(V^{\left[n\right]}\otimes i^{\left[n\right]}b^{\left[n\right]}\right) & = & \left(i^{\left[n\right]}b^{\left[n\right]}\otimes V^{\left[n\right]}\right)c_{V^{\left[n\right]},P^{\left[n\right]}},\label{form:bracket2}\end{eqnarray}
and\[
U^{\left[n+1\right]}=\frac{U^{\left[n\right]}}{\left(\left[\mathrm{Id}-i^{\left[n\right]}b^{\left[n\right]}\right]\left[P^{\left[n\right]}\right]\right)}.\]
When there is matter of confusion, we will use the notations $P_{V}^{\left[n\right]}:=P^{\left[n\right]},U_{V}^{\left[n\right]}:=U^{\left[n\right]},b_{V}^{\left[n\right]}:=b^{\left[n\right]},i_{V}^{\left[n\right]}:=i^{\left[n\right]}$
and $i_{U}^{V}:=i_{U}$ as well. 

We will denote by $\pi_{n}^{n+1}:U^{\left[n\right]}\rightarrow U^{\left[n+1\right]}$
the canonical projection for each $n\in\mathbb{N}$.
\end{defn}

Next we give the definition of bracket for a braided vector
space.
\begin{defn}
A \textbf{braided bracket }or simply a\textbf{ bracket} for a braided
vector space $\left(V,c\right)$ is a family $b:=\left(b^{\left[n\right]}\right)_{n\in\mathbb{N}}$
of maps such that, for each $n\in\mathbb{N}$, the family $\left(b^{\left[i\right]}\right)_{0\leq i\leq n}$
is a bracket of rank $n$ for $\left(V,c\right)$. This way we get
a direct system of $\left(V,c\right)$-bialgebras \[
\left(U^{\left[0\right]},i^{\left[0\right]}\right)\overset{\pi_{0}^{1}}{\rightarrow}\left(U^{\left[1\right]},i^{\left[1\right]}\right)\overset{\pi_{1}^{2}}{\rightarrow}\left(U^{\left[2\right]},i^{\left[2\right]}\right)\overset{\pi_{2}^{3}}{\rightarrow}\cdots.\]
 The direct limit of the direct system of braided bialgebras \[
U^{\left[0\right]}\overset{\pi_{0}^{1}}{\rightarrow}U^{\left[1\right]}\overset{\pi_{1}^{2}}{\rightarrow}U^{\left[2\right]}\overset{\pi_{2}^{3}}{\rightarrow}\cdots.\]
 will be denoted by \[
\left(U\left(V,c,b\right),\pi_{U}\right):=\left(U^{\left[\infty\right]},\pi_{n}^{\infty}\right):=\underrightarrow{\lim}U\left(V,c,\left(b^{\left[i\right]}\right)_{0\leq i\leq n-1}\right).\]
 By Lemma \ref{lem: lim of br bialg}, $U\left(V,c,b\right)$ carries
a natural braided bialgebra structure that makes it the direct limit
of $\left(\left(U^{\left[u\right]}\right)_{u\in\mathbb{N}},\left(\pi_{u}^{v}\right)_{u,v\in\mathbb{N}}\right)$
as a direct system of braided bialgebras. $U\left(V,c,b\right)$ will
be called \textbf{the} \textbf{universal enveloping algebra of} $\left(V,c,b\right).$

Denote by \[
i_{U}:V\rightarrow U\left(V,c,b\right)\]
 the (not necessarily injective) canonical map and set $P^{\left[\infty\right]}:=P\left(U^{\left[\infty\right]}\right).$
Note that $\mathrm{Im}\left(i_{U}\right)\subseteq P^{\left[\infty\right]}$
so that $i_{U}$ induces a morphism of braided vector spaces $i^{\left[\infty\right]}:V\rightarrow P^{\left[\infty\right]}$.

Clearly $\left(U^{\left[\infty\right]},i^{\left[\infty\right]}\right)$
is a $\left(V,c\right)$-bialgebra.

When there is matter of confusion, we will use the notation $i_{U}^{V}:=i_{U}$
as well. 
\end{defn}
Next aim is to describe the image of the map $i_{U}:V\rightarrow U\left(V,c,b\right)$. 
\begin{thm}
\label{teo: Im i_U}Let $\left(V,c\right)$ be a braided vector space
endowed with a bracket $b$. Then $P^{\left[\infty\right]}=\underrightarrow{\lim}P^{\left[n\right]}.$
Moreover $i^{\left[\infty\right]}:V\rightarrow P^{\left[\infty\right]}$
is surjective i.e. \[
\mathrm{Im}\left(i_{U}\right)=P\left(U\left(V,c,b\right)\right).\]
\end{thm}
\begin{proof}
The braided bialgebra homomorphism $\pi_{n-1}^{n}:=U^{\left[n-1\right]}\rightarrow U^{\left[n\right]}$
induces a direct system \[
P^{\left[0\right]}\overset{P\left(\pi_{0}^{1}\right)}{\rightarrow}P^{\left[1\right]}\overset{P\left(\pi_{1}^{2}\right)}{\rightarrow}P^{\left[2\right]}\overset{P\left(\pi_{2}^{3}\right)}{\rightarrow}\cdots.\]
 Let us prove that $P^{\left[\infty\right]}=\underrightarrow{\lim}P^{\left[n\right]}.$
Let $\delta^{\left[n\right]}:=\delta_{U^{\left[n\right]}}:U^{\left[n\right]}\rightarrow U^{\left[n\right]}\otimes U^{\left[n\right]}$
be defined as in \ref{nn: connected}. Then\begin{eqnarray*}
\left(\pi_{n-1}^{n}\otimes\pi_{n-1}^{n}\right)\delta^{\left[n-1\right]}\left(x\right) & = & \left(\pi_{n-1}^{n}\otimes\pi_{n-1}^{n}\right)\left(x\otimes1_{U^{\left[n-1\right]}}+1_{U^{\left[n-1\right]}}\otimes x-\Delta_{U^{\left[n-1\right]}}\left(x\right)\right)\\
 & = & \pi_{n-1}^{n}\left(x\right)\otimes1_{U^{\left[n\right]}}+1_{U^{\left[n\right]}}\otimes\pi_{n-1}^{n}\left(x\right)-\left(\pi_{n-1}^{n}\otimes\pi_{n-1}^{n}\right)\Delta_{U^{\left[n-1\right]}}\left(x\right)\\
 & = & \pi_{n-1}^{n}\left(x\right)\otimes1_{U^{\left[n\right]}}+1_{U^{\left[n\right]}}\otimes\pi_{n-1}^{n}\left(x\right)-\Delta_{U^{\left[n\right]}}\pi_{n-1}^{n}\left(x\right)\\
 & = & \delta^{\left[n\right]}\pi_{n-1}^{n}\left(x\right)\end{eqnarray*}
 so that $\left(\pi_{n-1}^{n}\otimes\pi_{n-1}^{n}\right)\delta^{\left[n-1\right]}=\delta^{\left[n\right]}\pi_{n-1}^{n}.$
Hence from the exact sequence\[
0\rightarrow P^{\left[n\right]}\rightarrow U^{\left[n\right]}\overset{\delta^{\left[n\right]}}{\longrightarrow}U^{\left[n\right]}\otimes U^{\left[n\right]}\]
 we get the exact sequence\[
0\rightarrow\underrightarrow{\lim}P^{\left[n\right]}\rightarrow\underrightarrow{\lim}U^{\left[n\right]}\overset{\underrightarrow{\lim}\delta^{\left[n\right]}}{\longrightarrow}\underrightarrow{\lim}\left(U^{\left[n\right]}\otimes U^{\left[n\right]}\right)\]
 which can be rewritten as\[
0\rightarrow\underrightarrow{\lim}P^{\left[n\right]}\rightarrow U^{\left[\infty\right]}\overset{\delta^{\left[\infty\right]}}{\longrightarrow}U^{\left[\infty\right]}\otimes U^{\left[\infty\right]}\]
 so that $P^{\left[\infty\right]}=\underrightarrow{\lim}P^{\left[n\right]}.$
For each $z\in P^{\left[\infty\right]},$ there is $t\in\mathbb{N}$
and $z_{t}\in P^{\left[t\right]}$ such that $z=\pi_{t}^{\infty}\left(z_{t}\right).$
Thus\[
z=\pi_{t}^{\infty}\left(z_{t}\right)=\pi_{t+1}^{\infty}\pi_{t}^{t+1}\left(z_{t}\right)\overset{z_{t}\in P^{\left[t\right]}}{=}\pi_{t+1}^{\infty}\pi_{t}^{t+1}i^{\left[t\right]}b^{\left[t\right]}\left(z_{t}\right)=i^{\left[\infty\right]}b^{\left[t\right]}\left(z_{t}\right)\in\mathrm{Im}\left(i^{\left[\infty\right]}\right)\]
 whence $P^{\left[\infty\right]}\subseteq\mathrm{Im}\left(i^{\left[\infty\right]}\right).$
Since $\mathrm{Im}\left(i^{\left[\infty\right]}\right)\subseteq P^{\left[\infty\right]}$
we conclude. 
\end{proof}
The map $i_{U}:V\rightarrow U\left(V,c,b\right)$ needs not to be
injective in general. Thus it makes sense to introduce the following
definition.
\begin{defn}
We say that $\left(V,c,b\right)$ is a \textbf{braided Lie algebra}
whenever 
\begin{itemize}
\item $\left(V,c\right)$ is a braided vector space; 
\item $b$ is a bracket on $\left(V,c\right);$ 
\item $i_{U}:V\rightarrow U\left(V,c,b\right)$ is injective (\textbf{implicit
Jacobi identity)}. 
\end{itemize}
\end{defn}

\begin{rem}
 The third condition above, was called implicit Jacobi identity having in mind \cite[Theorem 4.3]{AMS-MM2}.
\end{rem}

The following result will be crucial in proving the main theorem.

\begin{cor}
\label{cor:P(U)}Let $\left(V,c,b\right)$ be a braided Lie algebra.
The canonical map $i_{U}:V\rightarrow U\left(V,c,b\right)$ induces
an isomorphism of braided vector spaces between $V$ and $P\left(U\left(V,c,b\right)\right)$
namely the map $i^{\left[\infty\right]}:V\rightarrow P^{\left[\infty\right]}$. \end{cor}
\begin{proof}
Injectivity of $i_{U}:V\rightarrow U\left(V,c,b\right)$ implies Injectivity
of $i^{\left[\infty\right]}:V\rightarrow P^{\left[\infty\right]}.$
The conclusion follows by Theorem \ref{teo: Im i_U}. \end{proof}

We now give a different characterization of the implicit Jacobi identity.

\begin{prop}
\label{pro: i_U split}Let $\left(V,c\right)$ is a braided vector
space endowed with a bracket $b.$ The following assertions are equivalent. 
\begin{enumerate}
\item [$\left(i\right)$] the implicit Jacobi identity holds i.e. $\left(V,c,b\right)$
is a braided Lie algebra. 
\item [$\left(ii\right)$] $b^{\left[n\right]}i^{\left[n\right]}=\mathrm{Id}_{V}$
for all $n\in\mathbb{N}$. 
\item [$\left(iii\right)$] $b^{\left[0\right]}i^{\left[0\right]}=\mathrm{Id}_{V}$
and $b^{\left[n+1\right]}P\left(\pi_{n}^{n+1}\right)=b^{\left[n\right]}$
for all $n\in\mathbb{N}$. 
\end{enumerate}
Here $P\left(\pi_{n}^{n+1}\right):P^{\left[n\right]}\rightarrow P^{\left[n+1\right]}$
denotes the map induced by $\pi_{n}^{n+1}$.\end{prop}
\begin{proof}
$\left(i\right)\Rightarrow\left(ii\right)$ We have \begin{eqnarray*}
 &  & i_{U}\left(\mathrm{Id}_{V}-b^{\left[n\right]}i^{\left[n\right]}\right)\left(V\right)\\
 & = & P\left(\pi_{n}^{\infty}\right)i^{\left[n\right]}\left(\mathrm{Id}_{V}-b^{\left[n\right]}i^{\left[n\right]}\right)\left(V\right)\\
 & = & \pi_{n}^{\infty}\left(\mathrm{Id}_{U^{\left[n\right]}}-i^{\left[n\right]}b^{\left[n\right]}\right)i^{\left[n\right]}\left(V\right)\subseteq\pi_{n+1}^{\infty}\pi_{n}^{n+1}\left(\mathrm{Id}_{P^{\left[n\right]}}-i^{\left[n\right]}b^{\left[n\right]}\right)\left[P^{\left[n\right]}\right]=0\end{eqnarray*}
 so that $i_{U}\left(\mathrm{Id}_{V}-b^{\left[n\right]}i^{\left[n\right]}\right)=0.$
Then injectivity of $i_{U}$ yields $b^{\left[n\right]}i^{\left[n\right]}=\mathrm{Id}_{V}.$

$\left(ii\right)\Rightarrow\left(i\right)$ Since $i^{\left[n\right]}:V\rightarrow P^{\left[n\right]}$
is injective for all $n\in\mathbb{N}$, we get that $i^{\left[\infty\right]}=\underrightarrow{\lim}i^{\left[n\right]}:V\rightarrow\underrightarrow{\lim}P^{\left[n\right]}=P^{\left[\infty\right]}$
is injective too. Since $i^{\left[\infty\right]}$ is obtained by
corestricting $i_{U}$ to $P^{\left[\infty\right]},$ we have that
$i_{U}$ is injective too.

$\left(i\right)\Rightarrow\left(iii\right)$ For all $x\in P^{\left[n\right]}$
we have\begin{eqnarray*}
i_{U}\left[b^{\left[n\right]}-b^{\left[n+1\right]}P\left(\pi_{n}^{n+1}\right)\right]\left(x\right) & = & \left[i_{U}b^{\left[n\right]}-i_{U}b^{\left[n+1\right]}P\left(\pi_{n}^{n+1}\right)\right]\left(x\right)\\
 & = & \pi_{n+1}^{\infty}\pi_{n}^{n+1}i^{\left[n\right]}b^{\left[n\right]}\left(x\right)-\pi_{n+2}^{\infty}\pi_{n+1}^{n+2}i^{\left[n+1\right]}b^{\left[n+1\right]}P\left(\pi_{n}^{n+1}\right)\left(x\right)\\
 & = & \pi_{n+1}^{\infty}\pi_{n}^{n+1}\left(x\right)-\pi_{n+2}^{\infty}\pi_{n+1}^{n+2}P\left(\pi_{n}^{n+1}\right)\left(x\right)\\
 & = & \pi_{n}^{\infty}\left(x\right)-\pi_{n+1}^{\infty}P\left(\pi_{n}^{n+1}\right)\left(x\right)\\
 & = & \pi_{n}^{\infty}\left(x\right)-\pi_{n+1}^{\infty}\pi_{n}^{n+1}\left(x\right)=0\end{eqnarray*}
 so that $i_{U}\left[b^{\left[n\right]}-b^{\left[n+1\right]}P\left(\pi_{n}^{n+1}\right)\right]=0.$
Injectivity of $i_{U}$ implies $b^{\left[n\right]}=b^{\left[n+1\right]}P\left(\pi_{n}^{n+1}\right)$.
By $\left(i\right)\Rightarrow\left(ii\right)$ we have that $b^{\left[0\right]}i^{\left[0\right]}=\mathrm{Id}_{V}$.

$\left(iii\right)\Rightarrow\left(ii\right)$ We prove it by induction
on $n\in\mathbb{N}$. For $n=0$ there is nothing to prove. By definition
of $i^{\left[n+1\right]}$ and by induction hypothesis we get 
\begin{equation*}                                                   b^{\left[n+1\right]}i^{\left[n+1\right]}=b^{\left[n+1\right]}P\left(\pi_{n}^{n+1}\right)i^{\left[n\right]}=b^{\left[n\right]}i^{\left[n\right]}=\mathrm{Id}_{V}.\end{equation*}\end{proof}

In view of the previous result, we are able to give an alternative description of the universal
enveloping algebra of positive rank. 

\begin{cor}
\label{coro: ker b}Let $\left(V,c,b\right)$ be a braided Lie algebra.
Then\[
U^{\left[n+1\right]}=\frac{U^{\left[n\right]}}{\left(\mathrm{Ker}\left(b^{\left[n\right]}\right)\right)}\text{, for each }n\in\mathbb{N}.\]
\end{cor}
\begin{proof}
By Proposition \ref{pro: i_U split}, we have that $b_{V}^{\left[n\right]}i_{V}^{\left[n\right]}=\mathrm{Id}_{V}$
for each $n\in\mathbb{N}.$ This condition entails that $\mathrm{Im}\left(\mathrm{Id}-i_{V}^{\left[n\right]}b_{V}^{\left[n\right]}\right)=\mathrm{Ker}\left(b^{\left[n\right]}\right).$ 
\end{proof}

Next aim is to show that the universal enveloping algebra introduced
in \cite{Ardizzoni-Universal} can be understood as a universal enveloping algebra of rank one in the sense of the
present paper.

\begin{rem}
\label{rem:Ardi}Let $\left(V,c\right)$ be a braided vector space.
Denote by $E\left(V,c\right)$ the space spanned by primitive elements
of $T\left(V,c\right)$ of homogeneous degree at least two. In \cite[Definition 3.2]{Ardizzoni-Universal},
a bracket for $\left(V,c\right)$ was defined as a map $\beta:E\left(V,c\right)\rightarrow V$
such that \[
c\left(\beta\otimes V\right)=\left(V\otimes\beta\right)c_{E\left(V,c\right),V},\qquad c\left(V\otimes\beta\right)=\left(\beta\otimes V\right)c_{V,E\left(V,c\right)}.\]
The corresponding universal enveloping algebra was defined in \cite[Definition 3.5]{Ardizzoni-Universal}
as \[
\mathbb{U}\left(V,c,\beta\right)=\frac{T\left(V,c\right)}{\left(\left(\mathrm{Id}-\beta\right)\left[E\left(V,c\right)\right]\right)}.\]
When the canonical map $i_{\mathbb{U}}:V\rightarrow\mathbb{U}\left(V,c,\beta\right)$
is injective, the tern $\left(V,c,\beta\right)$ was called a braided
Lie algebra in \cite[Definition 4.1]{Ardizzoni-Universal}. Now the
primitive part of $\mathbb{U}\left(V,c,\beta\right)$ identifies with
$\left(V,c\right)$ itself if $\left(V,c\right)$ is in the class
$\mathcal{S}$ consisting of those braided vector spaces such that the kernel
of the canonical projection $p_{0}^{\infty}:T\left(V,c\right)\rightarrow\mathcal{B}\left(V,c\right)$, onto the Nichols algebra $\mathcal{B}\left(V,c\right)$, is the two-sided ideal of $T\left(V,c\right)$ generated by $E\left(V,c\right)$,
see \cite[Corollary 5.5]{Ardizzoni-Universal}. This property implies
that all primitively generated braided bialgebras whose infinitesimal
part is in $\mathcal{S}$ are universal enveloping algebras in the
above sense of their primitive part, see \cite[Theorem 5.7]{Ardizzoni-Universal}.
Significant braided vector spaces belong to the class $\mathcal{S}$,
see \cite{Ardizzoni-Sdeg}, where elements in $\mathcal{S}$ are called
braided vector spaces of strongness degree at most one. Anyway, the
fact that not all braided vector spaces are in $\mathcal{S}$, led
us to introduce in the present paper a different notion of braided
Lie algebra and universal enveloping algebra. Next result shows that
the old definition can be understood as a first step of the new one.\end{rem}
\begin{thm}
\label{teo:Ardi}Let $\left(V,c,b\right)$ be a braided Lie algebra.
Let $\beta:E\left(V,c\right)\rightarrow V$ be the restriction of
$b^{\left[0\right]}$ to $E\left(V,c\right)$. Then $\left(V,c,\beta\right)$
is a braided Lie algebra in the sense of Remark \ref{rem:Ardi} and
$\mathbb{U}\left(V,c,\beta\right)=U\left(V,c,b^{\left[0\right]}\right).$
Moreover, if the braided vector space $\left(V,c\right)$ is in the
class $\mathcal{S}$, then $U\left(V,c,b\right)=U\left(V,c,b^{\left[0\right]}\right)$
.\end{thm}
\begin{proof}
Since $i^{\left[0\right]}:V\rightarrow P^{\left[0\right]}$ is injective,
$V$ identifies with $V^{\left[0\right]}=i^{\left[0\right]}\left(V\right)$
so that (\ref{form:bracket1}) and (\ref{form:bracket2}) in case
$n=0$ rewrites as\[
c\left(b^{\left[0\right]}\otimes V\right)=\left(V\otimes b^{\left[0\right]}\right)c_{P^{\left[0\right]},V},\qquad c\left(V\otimes b^{\left[0\right]}\right)=\left(b^{\left[0\right]}\otimes V\right)c_{V,P^{\left[0\right]}}.\]
 Restricting the former equality to $E\left(V,c\right)\otimes V$
and the latter to $V\otimes E\left(V,c\right)$ yields\[
c\left(\beta\otimes V\right)=\left(V\otimes\beta\right)c_{E\left(V,c\right),V},\qquad c\left(V\otimes\beta\right)=\left(\beta\otimes V\right)c_{V,E\left(V,c\right)}\]
which means that $\beta$ is a bracket in the sense of Remark \ref{rem:Ardi}.
By \cite[Remark 2.9]{Ardizzoni-Universal}, $P^{\left[0\right]}=P\left(U^{\left[0\right]}\right)=V^{\left[0\right]}\oplus E\left(V,c\right).$
By Proposition \ref{pro: i_U split}, $b^{\left[n\right]}\circ i^{\left[n\right]}=\mathrm{Id}_{V}$
for all $n\in\mathbb{N}$ so that we have\begin{eqnarray*}
\left(\mathrm{Id}-i^{\left[0\right]}b^{\left[0\right]}\right)\left[P^{\left[0\right]}\right] & = & \left(\mathrm{Id}-i^{\left[0\right]}b^{\left[0\right]}\right)\left[V^{\left[0\right]}\oplus E\left(V,c\right)\right]\\
 & = & \left(\mathrm{Id}-i^{\left[0\right]}b^{\left[0\right]}\right)\left[V^{\left[0\right]}\right]+\left(\mathrm{Id}-i^{\left[0\right]}b^{\left[0\right]}\right)\left[E\left(V,c\right)\right]\\
 & = & \left(\mathrm{Id}-i^{\left[0\right]}b^{\left[0\right]}\right)i^{\left[0\right]}\left(V\right)+\left(\mathrm{Id}-i^{\left[0\right]}\beta\right)\left[E\left(V,c\right)\right]\\
 & = & \left(\mathrm{Id}-\beta\right)\left[E\left(V,c\right)\right]\end{eqnarray*}
 whence \[
U\left(V,c,b^{\left[0\right]}\right)=\frac{U^{\left[0\right]}}{\left(\left(\mathrm{Id}-i^{\left[0\right]}b^{\left[0\right]}\right)\left[P^{\left[0\right]}\right]\right)}=\frac{T\left(V,c\right)}{\left(\left(\mathrm{Id}-\beta\right)\left[E\left(V,c\right)\right]\right)}=\mathbb{U}\left(V,c,\beta\right).\]
Now, the canonical map $i_{\mathbb{U}}:V\rightarrow\mathbb{U}\left(V,c,\beta\right)$
corestricts to the map $i^{\left[1\right]}:V\rightarrow P^{\left[1\right]}$
which is split injective by $b^{\left[1\right]}$. Thus $i_{\mathbb{U}}$
is injective that is $\left(V,c,\beta\right)$ is a braided Lie algebra
in the sense of Remark \ref{rem:Ardi}. 

Assume now the braided vector space $\left(V,c\right)$ is in $\mathcal{S}$.
Then the space of primitive elements in $\mathbb{U}\left(V,c,\beta\right)$
identifies with $V$ via $i_{\mathbb{U}}$. Therefore, the projection
$\pi_{1}^{\infty}:U\left(V,c,b^{\left[0\right]}\right)\rightarrow U\left(V,c,b\right)$
is injective on primitive elements of $U\left(V,c,b^{\left[0\right]}\right)$.
By \cite[Lemma 5.3.3]{Montgomery} $\pi_{1}^{\infty}$ is injective
whence bijective.\end{proof}
\begin{lem}
\label{lem:UnivRank}Let $\left(V,c,b\right)$ be a braided Lie algebra
such that $U\left(V,c,b\right)=U\left(V,c,\left(b^{\left[t\right]}\right)_{0\leq t\leq n}\right)$
for some $n\in\mathbb{N}$. Then, for all $t\geq n+1$, $b^{\left[t\right]}=b^{\left[n+1\right]}$.
Moreover $b^{\left[n+1\right]}$ is bijective and it is the inverse
of $i^{\left[n+1\right]}:V\rightarrow P^{\left[n+1\right]}$.\end{lem}
\begin{proof}
Let $t\geq n+1$. By assumption, the canonical map $\pi_{n+1}^{\infty}:U^{\left[n+1\right]}\rightarrow U^{\left[\infty\right]}$
is the identity so that $\pi_{t}^{\infty}:U^{\left[t\right]}\rightarrow U^{\left[\infty\right]}$
is the identity too. Now, by Corollary \ref{cor:P(U)}, we know that
$i^{\left[\infty\right]}$ is bijective. By construction, $i^{\left[\infty\right]}=P\left(\pi_{t}^{\infty}\right)i^{\left[t\right]}=i^{\left[t\right]}$.
In particular $i^{\left[t\right]}$ is bijective and $i^{\left[t\right]}=i^{\left[n+1\right]}$.
From Proposition \ref{pro: i_U split}, we know that $b^{\left[t\right]}i^{\left[t\right]}=\mathrm{Id}_{V}$.
Thus $b^{\left[t\right]}$ is the inverse of $i^{\left[n+1\right]}:V\rightarrow P^{\left[n+1\right]}$
which is $b^{\left[n+1\right]}.$
\end{proof}

\section{The main result}\label{sec:mainresult}

In this section we associates a braided Lie algebra structure to the
infinitesimal part of a braided bialgebra. The braided bialgebra will
then be recovered as the universal enveloping algebra of its infinitesimal
braided Lie algebra.
\begin{prop}
\label{pro: alg is Lie}Any braided algebra $\left(A,c\right)$ admits
a canonical bracket $b=\left(b^{\left[n\right]}\right)_{n\in\mathbb{N}}$
such that $\left(A,c,b\right)$ is a braided Lie algebra.\end{prop}
\begin{proof}
Let us define recursively maps $b^{\left[n\right]}:P^{\left[n\right]}\rightarrow A$
such that $b^{\left[n\right]}i^{\left[n\right]}=\mathrm{Id}_{A}$
and (\ref{form:bracket1}) and (\ref{form:bracket2}) hold for $V=A$.

$n=0)$ By the universal property of $U^{\left[0\right]}=T\left(A,c\right)$
(see \cite[Theorem 1.17]{AMS-MM2}), $\mathrm{Id}_{A}$ can be lifted
to a unique braided algebra homomorphism $\varphi^{\left[0\right]}:U^{\left[0\right]}\rightarrow A$.
Let $b^{\left[0\right]}:P^{\left[0\right]}\rightarrow A$ be the restriction
of $\varphi^{\left[0\right]}$ to $P^{\left[0\right]}$. Clearly $b^{\left[0\right]}\circ i^{\left[0\right]}=\mathrm{Id}_{A}.$
Let us check that\begin{equation}
c_{U^{\left[0\right]},A^{\left[0\right]}}\left(i^{\left[0\right]}\varphi^{\left[0\right]}\otimes A^{\left[0\right]}\right)=\left(A^{\left[0\right]}\otimes i^{\left[0\right]}\varphi^{\left[0\right]}\right)c_{U^{\left[0\right]},A^{\left[0\right]}}\label{form:phi0}\end{equation}
holds (here $A^{\left[0\right]}:=i^{\left[0\right]}\left(A\right)$).
Note that the notation $c_{U^{\left[0\right]},A^{\left[0\right]}}$
makes sense as, by Proposition \ref{pro: cat-braid}, $A^{\left[0\right]}$
is a categorical subspace of $\left(U^{\left[0\right]},c_{U^{\left[0\right]}}\right).$
Let $j_{t}:A^{\otimes t}\rightarrow U^{\left[0\right]}$ be the canonical
injection. Then\begin{eqnarray*}
 &  & c_{U^{\left[0\right]},A^{\left[0\right]}}\left(i^{\left[0\right]}\varphi^{\left[0\right]}\otimes A^{\left[0\right]}\right)\left(j_{t}\otimes j_{0}\right)\\
 & = & c_{U^{\left[0\right]},A^{\left[0\right]}}\left(i^{\left[0\right]}\varphi^{\left[0\right]}j_{t}\otimes j_{0}\right)\\
 & = & c_{U^{\left[0\right]},A^{\left[0\right]}}\left(i^{\left[0\right]}\varphi^{\left[0\right]}j_{t}\otimes i^{\left[0\right]}\right)\\
 & = & \left(i^{\left[0\right]}\otimes i^{\left[0\right]}\right)c_{A}\left(\varphi^{\left[0\right]}j_{t}\otimes A\right)\\
 & = & \left(i^{\left[0\right]}\otimes i^{\left[0\right]}\right)c_{A}\left(m_{A}^{t-1}\otimes A\right)\\
 & \overset{\left(\ref{Br2}\right)}{=} & \left(i^{\left[0\right]}\otimes i^{\left[0\right]}\right)\left(A\otimes m_{A}^{t-1}\right)c_{A^{\otimes t},A}\\
 & = & \left(i^{\left[0\right]}\otimes i^{\left[0\right]}\right)\left(A\otimes\varphi^{\left[0\right]}j_{t}\right)c_{A^{\otimes t},A}\\
 & = & \left(A^{\left[0\right]}\otimes i^{\left[0\right]}\varphi^{\left[0\right]}\right)\left(j_{0}\otimes j_{t}\right)c_{A^{\otimes t},A}\\
 & = & \left(A^{\left[0\right]}\otimes i^{\left[0\right]}\varphi^{\left[0\right]}\right)c_{U^{\left[0\right]},U^{\left[0\right]}}\left(j_{t}\otimes j_{0}\right)\\
 & = & \left(A^{\left[0\right]}\otimes i^{\left[0\right]}\varphi^{\left[0\right]}\right)c_{U^{\left[0\right]},A^{\left[0\right]}}\left(j_{t}\otimes j_{0}\right).\end{eqnarray*}
 Since this equality holds for all $t\in\mathbb{N},$ and $A^{\left[0\right]}=\mathrm{Im}\left(i^{\left[0\right]}\right)=\mathrm{Im}\left(j_{0}\right),$
we get (\ref{form:phi0}). Now we can prove that (\ref{form:bracket1})
is fulfilled for $n=0.$ In fact, for each $x\in P^{\left[0\right]}$
and $a\in A^{\left[0\right]}$ we have\begin{eqnarray*}
 &  & c_{A^{\left[0\right]}}\left(i^{\left[0\right]}b^{\left[0\right]}\otimes A^{\left[0\right]}\right)\left(x\otimes a\right)\\
 & = & c_{U^{\left[0\right]},A^{\left[0\right]}}\left(i^{\left[0\right]}\varphi^{\left[0\right]}\otimes A^{\left[0\right]}\right)\left(x\otimes a\right)\\
 & \overset{\eqref{form:phi0}}{=} & \left(A^{\left[0\right]}\otimes i^{\left[0\right]}\varphi^{\left[0\right]}\right)c_{U^{\left[0\right]},A^{\left[0\right]}}\left(x\otimes a\right)\\
 & = & \left(A^{\left[0\right]}\otimes i^{\left[0\right]}b^{\left[0\right]}\right)c_{P^{\left[0\right]},A^{\left[0\right]}}\left(x\otimes a\right).\end{eqnarray*}
 Similarly one proves (\ref{form:bracket2}).

$n\geq1)$ We have\begin{eqnarray*}
 &  & \varphi^{\left[n-1\right]}\left(\mathrm{Id}-i^{\left[n-1\right]}b^{\left[n-1\right]}\right)\left[P^{\left[n-1\right]}\right]\\
 & = & \left(\varphi^{\left[n-1\right]}-\varphi^{\left[n-1\right]}i^{\left[n-1\right]}b^{\left[n-1\right]}\right)\left[P^{\left[n-1\right]}\right]\\
 & = & \left(\varphi^{\left[n-1\right]}-b^{\left[n-1\right]}\right)\left[P^{\left[n-1\right]}\right]=\left(\varphi^{\left[n-1\right]}-\varphi^{\left[n-1\right]}\right)\left[P^{\left[n-1\right]}\right]=0\end{eqnarray*}
so that $\varphi^{\left[n-1\right]}:U^{\left[n-1\right]}\rightarrow A$
quotients to a map $\varphi^{\left[n\right]}:U^{\left[n\right]}\rightarrow A$.
Let $b^{\left[n\right]}:P^{\left[n\right]}\rightarrow A$ be the restriction
of $\varphi^{\left[n\right]}$ to $P^{\left[n\right]}$. Clearly $b^{\left[n\right]}i^{\left[n\right]}=\mathrm{Id}_{A}$.
Let us prove that (\ref{form:bracket1}) holds. Let $x\in P^{\left[n\right]}$
and $y\in A^{\left[n\right]}.$ Then there is $x_{A}\in U^{\left[0\right]}$
such that $x=\pi_{0}^{n}\left(x_{A}\right)$ and $y_{A}\in A$ such
that $y=i^{\left[n\right]}\left(y_{A}\right)=\pi_{0}^{n}i^{\left[0\right]}\left(y_{A}\right)$.
Whence\begin{eqnarray*}
 &  & c_{A^{\left[n\right]}}\left(i^{\left[n\right]}b^{\left[n\right]}\otimes A^{\left[n\right]}\right)\left(x\otimes y\right)\\
 & = & c_{U^{\left[n\right]}}\left(i^{\left[n\right]}b^{\left[n\right]}\otimes A^{\left[n\right]}\right)\left(\pi_{0}^{n}\left(x_{A}\right)\otimes\pi_{0}^{n}i^{\left[0\right]}\left(y_{A}\right)\right)\\
 & = & c_{U^{\left[n\right]}}\left(i^{\left[n\right]}\varphi^{\left[n\right]}\pi_{0}^{n}\left(x_{A}\right)\otimes\pi_{0}^{n}i^{\left[0\right]}\left(y_{A}\right)\right)\\
 & = & c_{U^{\left[n\right]}}\left(\pi_{0}^{n}i^{\left[0\right]}\varphi^{\left[0\right]}\left(x_{A}\right)\otimes\pi_{0}^{n}i^{\left[0\right]}\left(y_{A}\right)\right)\\
 & = & \left(\pi_{0}^{n}\otimes\pi_{0}^{n}\right)c_{U^{\left[0\right]}}\left(i^{\left[0\right]}\varphi^{\left[0\right]}\otimes U^{\left[0\right]}\right)\left(x_{A}\otimes i^{\left[0\right]}\left(y_{A}\right)\right)\\
 & = & \left(\pi_{0}^{n}\otimes\pi_{0}^{n}\right)c_{U^{\left[0\right]},A^{\left[0\right]}}\left(i^{\left[0\right]}\varphi^{\left[0\right]}\otimes A^{\left[0\right]}\right)\left(x_{A}\otimes i^{\left[0\right]}\left(y_{A}\right)\right)\\
 & \overset{\eqref{form:phi0}}{=} & \left(\pi_{0}^{n}\otimes\pi_{0}^{n}\right)\left(A^{\left[0\right]}\otimes i^{\left[0\right]}\varphi^{\left[0\right]}\right)c_{U^{\left[0\right]},A^{\left[0\right]}}\left(x_{A}\otimes i^{\left[0\right]}\left(y_{A}\right)\right)\\
 & = & \left(\pi_{0}^{n}\otimes\pi_{0}^{n}i^{\left[0\right]}\varphi^{\left[0\right]}\right)c_{U^{\left[0\right]},A^{\left[0\right]}}\left(x_{A}\otimes i^{\left[0\right]}\left(y_{A}\right)\right)\\
 & = & \left(\pi_{0}^{n}\otimes i^{\left[n\right]}\varphi^{\left[0\right]}\right)c_{U^{\left[0\right]}}\left(x_{A}\otimes i^{\left[0\right]}\left(y_{A}\right)\right)\\
 & = & \left(\pi_{0}^{n}\otimes i^{\left[n\right]}\varphi^{\left[n\right]}\pi_{0}^{n}\right)c_{U^{\left[0\right]}}\left(x_{A}\otimes i^{\left[0\right]}\left(y_{A}\right)\right)\\
 & = & \left(A^{\left[n\right]}\otimes i^{\left[n\right]}b^{\left[n\right]}\right)c_{U^{\left[n\right]}}\left(\pi_{0}^{n}\left(x_{A}\right)\otimes\pi_{0}^{n}i^{\left[0\right]}\left(y_{A}\right)\right)\\
 & = & \left(A^{\left[n\right]}\otimes i^{\left[n\right]}b^{\left[n\right]}\right)c_{P^{\left[n\right]},A^{\left[n\right]}}\left(x\otimes y\right)\end{eqnarray*}
 Similarly one proves (\ref{form:bracket2}). We have so proved that
$b:=\left(b^{\left[n\right]}\right)_{n\in\mathbb{N}}$ is a bracket
for $\left(A,c\right).$ Since, $b^{\left[n\right]}i^{\left[n\right]}=\mathrm{Id}_{A}$
for each $n\in\mathbb{N}$, by Proposition \ref{pro: i_U split}, we conclude
that $\left(A,c,b\right)$ is a braided Lie algebra. \end{proof}
\begin{claim}
\label{claim: P_h}Let $h:\left(V,c_{V}\right)\rightarrow\left(W,c_{W}\right)$
be a morphism of braided vector spaces. Assume that $\left(V,c_{V}\right)$
and $\left(W,c_{W}\right)$ are endowed with a bracket $b_{V}$ and
$b_{W}$ respectively. By the universal property of the tensor algebra,
$h$ induces a braided bialgebra homomorphism \[
U_{h}^{\left[0\right]}:U_{V}^{\left[0\right]}\rightarrow U_{W}^{\left[0\right]}\]
such that $U_{h}^{\left[0\right]}i_{V}^{\left[0\right]}=i_{W}^{\left[0\right]}h$
where $i_{W}^{\left[0\right]}$ is thought as a map into $U_{W}^{\left[0\right]}$
and $i_{V}^{\left[0\right]}$ as a map into $U_{V}^{\left[0\right]}$.
The map $U_{h}^{\left[0\right]}$ preserves primitives so that it
gives rise to a map \[
P_{h}^{\left[0\right]}:P_{V}^{\left[0\right]}\rightarrow P_{W}^{\left[0\right]}\]
 such that $P_{h}^{\left[0\right]}i_{V}^{\left[0\right]}=i_{W}^{\left[0\right]}h.$
Let $\pi_{0}^{1}:U_{V}^{\left[0\right]}\rightarrow U_{V}^{\left[1\right]}$
and $\Pi_{0}^{1}:U_{W}^{\left[0\right]}\rightarrow U_{W}^{\left[1\right]}$
be the canonical projections. Suppose \begin{equation}
hb_{V}^{\left[n\right]}=b_{W}^{\left[n\right]}P_{h}^{\left[n\right]}\label{form:morphbrackets}\end{equation}
 holds for $n=0.$ Since\begin{eqnarray*}
\Pi_{0}^{1}U_{h}^{\left[0\right]}\left(\mathrm{Id}-i_{V}^{\left[0\right]}b_{V}^{\left[0\right]}\right)\left[P_{V}^{\left[0\right]}\right] & = & \Pi_{0}^{1}\left(U_{h}^{\left[0\right]}-U_{h}^{\left[0\right]}i_{V}^{\left[0\right]}b_{V}^{\left[0\right]}\right)\left[P_{V}^{\left[0\right]}\right]\\
 & = & \Pi_{0}^{1}\left(P_{h}^{\left[0\right]}-i_{W}^{\left[0\right]}hb_{V}^{\left[0\right]}\right)\left[P_{V}^{\left[0\right]}\right]\\
 & = & \Pi_{0}^{1}\left(P_{h}^{\left[0\right]}-i_{W}^{\left[0\right]}b_{W}^{\left[0\right]}P_{h}^{\left[0\right]}\right)\left[P_{V}^{\left[0\right]}\right]\\
 & = & \Pi_{0}^{1}\left(\mathrm{Id}-i_{W}^{\left[0\right]}b_{W}^{\left[0\right]}\right)P_{h}^{\left[0\right]}\left[P_{V}^{\left[0\right]}\right]\subseteq\Pi_{0}^{1}\left(\mathrm{Id}-i_{W}^{\left[0\right]}b_{W}^{\left[0\right]}\right)\left[P_{W}^{\left[0\right]}\right]=0\end{eqnarray*}
we have that $U_{h}^{\left[0\right]}$ induces a map $U_{h}^{\left[1\right]}:U_{V}^{\left[1\right]}\rightarrow U_{W}^{\left[1\right]}$
such that $U_{h}^{\left[1\right]}\pi_{0}^{1}=\Pi_{0}^{1}U_{h}^{\left[0\right]}.$
Using that $\Pi_{0}^{1}$ is a surjective braided bialgebra homomorphism,
one can check that $U_{h}^{\left[1\right]}$ is a braided bialgebra
homomorphism too. Note that \[
U_{h}^{\left[1\right]}i_{V}^{\left[1\right]}=U_{h}^{\left[1\right]}\pi_{0}^{1}i_{V}^{\left[0\right]}=\Pi_{0}^{1}U_{h}^{\left[0\right]}i_{V}^{\left[0\right]}=\Pi_{0}^{1}i_{W}^{\left[0\right]}h=i_{W}^{\left[1\right]}h.\]
 As in the case of $U_{h}^{\left[0\right]}$, the map $U_{h}^{\left[1\right]}$
preserves primitives so that it gives rise to a map \[
P_{h}^{\left[1\right]}:P_{V}^{\left[1\right]}\rightarrow P_{W}^{\left[1\right]}.\]
 Now supposing \eqref{form:morphbrackets} holds for $n=1,$ one can
define $P_{h}^{\left[2\right]}$ and so on. \end{claim}
\begin{defn}
A \textbf{morphism of braided brackets }or simply\textbf{ a morphism
of brackets} $h:\left(V,c_{V},b_{V}\right)\rightarrow\left(W,c_{W},b_{W}\right)$
is a morphism of braided vector spaces $h:\left(V,c_{V}\right)\rightarrow\left(W,c_{W}\right)$
such that (\ref{form:morphbrackets}) holds for every $n\in\mathbb{N}$,
where $P_{h}^{\left[n\right]}$ is constructed as in \ref{claim: P_h}.
Note that $P_{h}^{\left[0\right]}$ always exists while, for $t\geq1,$
$P_{h}^{\left[t\right]}$ exists as (\ref{form:morphbrackets}) holds
for $n=t-1.$ \end{defn}
\begin{lem}
\label{lem: subLie}Let $\left(W,c_{W},b_{W}\right)$ a braided Lie
algebra and let $h:\left(V,c_{V}\right)\rightarrow\left(W,c_{W}\right)$
be an injective morphism of braided vector spaces. If, for each $n\in\mathbb{N}$,
there exists a $K$-linear map $b_{V}^{\left[n\right]}:P_{V}^{\left[n\right]}\rightarrow V$
such that (\ref{form:morphbrackets}) holds, where $P_{h}^{\left[n\right]}:P_{V}^{\left[n\right]}\rightarrow P_{W}^{\left[n\right]}$
is constructed step by step as in \ref{claim: P_h}, then $\left(V,c_{V},b_{V}\right)$
is a braided Lie algebra too and $h$ is a morphism of braided brackets. \end{lem}
\begin{proof}
Construct $U_{h}^{\left[0\right]}:U_{V}^{\left[0\right]}\rightarrow U_{W}^{\left[0\right]}$
and $P_{h}^{\left[0\right]}:P_{V}^{\left[0\right]}\rightarrow P_{W}^{\left[0\right]}$
as in \ref{claim: P_h}. Note that $U_{h}^{\left[0\right]}$ restricts
to a map $h^{\left[0\right]}:V^{\left[0\right]}\rightarrow W^{\left[0\right]}$
such that $h^{\left[0\right]}i_{V}^{\left[0\right]}=i_{W}^{\left[0\right]}h$
where $i_{V}^{\left[0\right]}$ is thought as a map into $V^{\left[0\right]}$
and $i_{W}^{\left[0\right]}$ as a map into $W^{\left[0\right]}$.
For all $w\in P_{V}^{\left[0\right]}$ and $y\in V^{\left[0\right]}$
we have\begin{eqnarray*}
 &  & c_{P_{W}^{\left[0\right]},W^{\left[0\right]}}\left(P_{h}^{\left[0\right]}\otimes h^{\left[0\right]}\right)\left(w\otimes y\right)\\
 & = & c_{U_{W}^{\left[0\right]},U_{W}^{\left[0\right]}}\left(U_{h}^{\left[0\right]}\otimes U_{h}^{\left[0\right]}\right)\left(w\otimes y\right)\\
 & = & \left(U_{h}^{\left[0\right]}\otimes U_{h}^{\left[0\right]}\right)c_{U_{V}^{\left[0\right]},U_{V}^{\left[0\right]}}\left(w\otimes y\right)=\left(h^{\left[0\right]}\otimes P_{h}^{\left[0\right]}\right)c_{P_{V}^{\left[0\right]},V^{\left[0\right]}}\left(w\otimes y\right).\end{eqnarray*}
Thus we have \begin{equation}
c_{P_{W}^{\left[0\right]},W^{\left[0\right]}}\left(P_{h}^{\left[0\right]}\otimes h^{\left[0\right]}\right)=\left(h^{\left[0\right]}\otimes P_{h}^{\left[0\right]}\right)c_{P_{V}^{\left[0\right]},V^{\left[0\right]}}.\label{form:h0}\end{equation}
 Hence, we have\begin{eqnarray*}
 &  & \left(h^{\left[0\right]}\otimes h^{\left[0\right]}\right)c_{V^{\left[n\right]}}\left(i_{V}^{\left[0\right]}b_{V}^{\left[0\right]}\otimes V^{\left[0\right]}\right)\\
 & = & c_{W^{\left[0\right]}}\left(h^{\left[0\right]}\otimes h^{\left[0\right]}\right)\left(i_{V}^{\left[0\right]}b_{V}^{\left[0\right]}\otimes V^{\left[0\right]}\right)\\
 & = & c_{W^{\left[0\right]}}\left(i_{W}^{\left[0\right]}hb_{V}^{\left[0\right]}\otimes h^{\left[0\right]}\right)\\
 & \overset{(\ref{form:morphbrackets})}{=} & c_{W^{\left[0\right]}}\left(i_{W}^{\left[0\right]}b_{W}^{\left[0\right]}P_{h}^{\left[0\right]}\otimes h^{\left[0\right]}\right)\\
 & = & c_{W^{\left[0\right]}}\left(i_{W}^{\left[0\right]}b_{W}^{\left[0\right]}\otimes W^{\left[0\right]}\right)\left(P_{h}^{\left[0\right]}\otimes h^{\left[0\right]}\right)\\
 & \overset{\left(\ref{form:bracket1}\right)}{=} & \left(W^{\left[0\right]}\otimes i_{W}^{\left[0\right]}b_{W}^{\left[0\right]}\right)c_{P_{W}^{\left[0\right]},W^{\left[0\right]}}\left(P_{h}^{\left[0\right]}\otimes h^{\left[0\right]}\right)\\
 & \overset{\left(\ref{form:h0}\right)}{=} & \left(h^{\left[0\right]}\otimes i_{W}^{\left[0\right]}b_{W}^{\left[0\right]}P_{h}^{\left[0\right]}\right)c_{P_{V}^{\left[0\right]},V^{\left[0\right]}}\\
 & \overset{(\ref{form:morphbrackets})}{=} & \left(h^{\left[0\right]}\otimes i_{W}^{\left[0\right]}hb_{V}^{\left[0\right]}\right)c_{P_{V}^{\left[0\right]},V^{\left[0\right]}}\\
 & = & \left(h^{\left[0\right]}\otimes h^{\left[0\right]}\right)\left(V^{\left[0\right]}\otimes i_{V}^{\left[0\right]}b_{V}^{\left[0\right]}\right)c_{P_{V}^{\left[0\right]},V^{\left[0\right]}}\end{eqnarray*}
 so that\[
\left(h^{\left[0\right]}\otimes h^{\left[0\right]}\right)c_{V^{\left[n\right]}}\left(i_{V}^{\left[0\right]}b_{V}^{\left[0\right]}\otimes V^{\left[0\right]}\right)=\left(h^{\left[0\right]}\otimes h^{\left[0\right]}\right)\left(V^{\left[0\right]}\otimes i_{V}^{\left[0\right]}b_{V}^{\left[0\right]}\right)c_{P_{V}^{\left[0\right]},V^{\left[0\right]}}.\]
Now, since $\left(W,c_{W},b_{W}\right)$ is a braided Lie algebra,
by Proposition \ref{pro: i_U split}, we have that $i_{W}^{\left[0\right]}$
is split injective. From this fact, injectivity of $h$ and equality
$h^{\left[0\right]}i_{V}^{\left[0\right]}=i_{W}^{\left[0\right]}h,$
one easily gets that $\mathrm{Ker}\left(h^{\left[0\right]}\right)=0$
so that $h^{\left[0\right]}$ is injective. Thus, from the last displayed
equality, we infer that $b_{V}^{\left[0\right]}$ fulfills (\ref{form:bracket1}).
Similarly (\ref{form:bracket2}) is satisfied. Hence $b_{V}^{\left[0\right]}$
is a bracket for $U_{V}^{\left[0\right]}$ and one can construct $U_{V}^{\left[1\right]}.$
Thus we can construct $U_{h}^{\left[1\right]}:U_{V}^{\left[1\right]}\rightarrow U_{W}^{\left[1\right]}$
and $P_{h}^{\left[1\right]}:P_{V}^{\left[1\right]}\rightarrow P_{W}^{\left[1\right]}$
as in \ref{claim: P_h} and $h^{\left[1\right]}:V^{\left[1\right]}\rightarrow W^{\left[1\right]}$
as we did for $h^{\left[0\right]}$ above. Then one can proceed as
in the case of $U_{h}^{\left[0\right]}$ to prove that $b_{V}^{\left[1\right]}$
fulfills (\ref{form:bracket1}) and (\ref{form:bracket2}). Going
on this way, we get that $b_{V}=\left(b_{V}^{\left[n\right]}\right)_{n\in\mathbb{N}}$
is a bracket on $\left(V,c_{V}\right).$ By (\ref{form:morphbrackets})$,$
$h$ is a morphism of braided brackets. Now, let \[
U_{h}^{\left[\infty\right]}:=\underrightarrow{\lim}U_{h}^{\left[n\right]}:U\left(V,c_{V},b_{V}\right)=U_{V}^{\left[\infty\right]}\rightarrow U_{W}^{\left[\infty\right]}=U\left(W,c_{W},b_{W}\right).\]
Denoting by $i_{U}^{V}:V\rightarrow U\left(V,c_{V},b_{V}\right)$
and $i_{U}^{W}:W\rightarrow U\left(W,c_{W},b_{W}\right)$ the canonical
maps, we get that $U_{h}^{\left[\infty\right]}i_{U}^{V}=i_{U}^{W}h.$
Since both $i_{U}^{W}$ and $h$ are injective, we deduce that $i_{U}^{V}$
is injective too. Therefore $\left(V,c_{V},b_{V}\right)$ is a braided
Lie algebra. \end{proof}
\begin{thm}
\label{thm:Infinitesimal}Let $A$ be a braided bialgebra with infinitesimal
part $(P,c_{P})$. Then the braided Lie algebra structure of $A$
of Proposition \ref{pro: alg is Lie} induces a braided Lie algebra
structure $\left(P,c_{P},b_{P}\right)$. \end{thm}
\begin{proof}
Let $\varphi^{\left[n\right]}:U_{A}^{\left[n\right]}\rightarrow A$
be defined as in the proof of Proposition \ref{pro: alg is Lie}.
In order to avoid ambiguities on the notation $P^{\left[n\right]},$
we set $V:=P$. By Proposition \ref{pro: alg is Lie}, $A$ has a
braided Lie algebra structure, say$\left(A,c_{A},b_{A}\right)$, and
by Proposition \ref{pro: i_U split}, $b_{A}^{\left[n\right]}i_{A}^{\left[n\right]}=\mathrm{Id}_{A}$
for all $n\in\mathbb{N}$. The canonical injection $\lambda:V\rightarrow A$
is a morphism of braided vector spaces. We will apply Lemma \ref{lem: subLie}
and use the notations in its proof. For any $n\in\mathbb{N}$, consider
the map \[
b_{A}^{\left[n\right]}P_{\lambda}^{\left[n\right]}:P_{V}^{\left[n\right]}\rightarrow A.\]
 By construction of $b_{A}^{\left[n\right]},$ for $x\in P_{V}^{\left[n\right]}$
we have\[
b_{A}^{\left[n\right]}P_{\lambda}^{\left[n\right]}\left(x\right)=\varphi^{\left[n\right]}P_{\lambda}^{\left[n\right]}\left(x\right)=\varphi^{\left[n\right]}U_{\lambda}^{\left[n\right]}\left(x\right).\]
 Since $\varphi^{\left[n\right]}U_{\lambda}^{\left[n\right]}$ is
a braided bialgebra homomorphism, it preserves primitive elements
so that $\varphi^{\left[n\right]}U_{\lambda}^{\left[n\right]}\left(x\right)\in V.$
Thus $b_{A}^{\left[n\right]}P_{\lambda}^{\left[n\right]}\left(x\right)\in V$
whence there exists a unique map $b_{V}^{\left[n\right]}:P_{V}^{\left[n\right]}\rightarrow V$
such that $\lambda b_{V}^{\left[n\right]}=b_{A}^{\left[n\right]}P_{\lambda}^{\left[n\right]}.$
By Lemma \ref{lem: subLie}, $\left(V,c_{V},b_{V}\right)$ is a braided
Lie algebra too.\end{proof}
\begin{defn}
With the same assumptions and notations as in Theorem \ref{thm:Infinitesimal},
$\left(P,c_{P},b_{P}\right)$ will be called the \textbf{infinitesimal
braided Lie algebra of} $A.$ \end{defn}

We are now able to state the universal property of the universal enveloping algebra.

\begin{thm}[\textbf{The universal property}]
\label{teo: univ U}Let $A$ be a braided bialgebra with infinitesimal
braided Lie algebra $\left(P,c_{P},b_{P}\right)$. Then every morphism
of braided brackets $h:\left(V,c,b\right)\rightarrow\left(P,c_{P},b_{P}\right)$
can be lifted to a morphism of braided bialgebras $\overline{h}:U\left(V,c,b\right)\rightarrow A.$ \end{thm}
\begin{proof}
Let $\varphi^{\left[n\right]}:U_{A}^{\left[n\right]}\rightarrow A$
be defined as in the proof of Proposition \ref{pro: alg is Lie}.
Let $\lambda:P\rightarrow A$ be the canonical map. Set $h^{\left[n\right]}:=\varphi^{\left[n\right]}U_{\lambda h}^{\left[n\right]}:U_{V}^{\left[n\right]}\rightarrow A.$
Let $\pi_{n-1}^{n}:U_{V}^{\left[n-1\right]}\rightarrow U_{V}^{\left[n\right]}$
and $\Pi_{n-1}^{n}:U_{A}^{\left[n-1\right]}\rightarrow U_{A}^{\left[n\right]}$
be the canonical projections. We have\[
h^{\left[n\right]}\pi_{n-1}^{n}=\varphi^{\left[n\right]}U_{\lambda h}^{\left[n\right]}\pi_{n-1}^{n}=\varphi^{\left[n\right]}\Pi_{n-1}^{n}U_{\lambda h}^{\left[n-1\right]}=\varphi^{\left[n-1\right]}U_{\lambda h}^{\left[n-1\right]}=h^{\left[n-1\right]}.\]
Since $U\left(V,c,b\right)$ is the direct limit of the direct system
\[
U_{V}^{\left[0\right]}\overset{\pi_{0}^{1}}{\rightarrow}U_{V}^{\left[1\right]}\overset{\pi_{1}^{2}}{\rightarrow}U_{V}^{\left[2\right]}\overset{\pi_{2}^{3}}{\rightarrow}\cdots\]
of braided bialgebras, we get that the compatible family $\left(h^{\left[n\right]}:U_{V}^{\left[n\right]}\rightarrow A\right)_{n\in\mathbb{N}}$
induces a braided bialgebra homomorphism $\overline{h}:U\left(V,c,b\right)\rightarrow A$
such that $\overline{h}\pi_{n}^{\infty}=h^{\left[n\right]}.$ In particular
\[
\overline{h}i_{U}=\overline{h}\pi_{0}^{\infty}i_{V}^{\left[0\right]}=h^{\left[0\right]}i_{V}^{\left[0\right]}=\varphi^{\left[0\right]}U_{\lambda h}^{\left[0\right]}i_{V}^{\left[0\right]}=\varphi^{\left[0\right]}i_{A}^{\left[0\right]}\lambda h=\lambda h.\] \end{proof}

The following Theorem \ref{teo:infpartU} and Theorem \ref{teo: generated} constitutes together a Milnor-Moore type theorem for primitively generated braided bialgebras, compare with \cite[Theorem 5.18]{Milnor-Moore}.

\begin{thm}
\label{teo:infpartU}Let $\left(V,c,b\right)$ be e braided Lie algebra.
Then the infinitesimal braided Lie algebra of $U\left(V,c,b\right)$
is $\left(V,c,b\right)$.\end{thm}
\begin{proof}
Let $U:=U\left(V,c,b\right)$ and let $\left(W,c_{W},b_{W}\right)$
be its infinitesimal braided Lie algebra. By Corollary \ref{cor:P(U)},
the canonical map $i_{U}:V\rightarrow U$ induces an isomorphism of
braided vector spaces between $V$ and $W$ namely the map $h:=i^{\left[\infty\right]}:V\rightarrow W$.
We only have to check this is a morphism of braided brackets. Let
us prove that \eqref{form:morphbrackets} holds by complete induction
on $n\in\mathbb{N}$. Let $\varphi^{\left[n\right]}:U_{U}^{\left[n\right]}\rightarrow U$
be defined as in the proof of Proposition \ref{pro: alg is Lie}.
For $n=0$ and $x\in P_{V}^{\left[0\right]}$, we have \begin{eqnarray*}
b_{W}^{\left[0\right]}P_{h}^{\left[0\right]}\left(x\right) & = & b_{U}^{\left[0\right]}P_{h}^{\left[0\right]}\left(x\right)=\varphi^{\left[0\right]}P_{h}^{\left[0\right]}\left(x\right)=\varphi^{\left[0\right]}U_{h}^{\left[0\right]}\left(x\right)=\pi_{0}^{\infty}\left(x\right)=\pi_{0}^{\infty}i^{\left[0\right]}b^{\left[0\right]}\left(x\right)=hb^{\left[0\right]}\left(x\right)\end{eqnarray*}
so that $b_{W}^{\left[0\right]}P_{h}^{\left[0\right]}=hb^{\left[0\right]}$.
Let $n>0$ and assume the statement true for all $0\leq i\leq n-1$.
Then, following \ref{claim: P_h}, we know that $P_{h}^{\left[n\right]}$
exists. For $x\in P_{V}^{\left[n\right]}$, we have \begin{eqnarray*}
b_{W}^{\left[n\right]}P_{h}^{\left[n\right]}\left(x\right) & = & b_{U}^{\left[n\right]}P_{h}^{\left[n\right]}\left(x\right)=\varphi^{\left[n\right]}P_{h}^{\left[n\right]}\left(x\right)=\varphi^{\left[n\right]}U_{h}^{\left[n\right]}\left(x\right)=\pi_{n}^{\infty}\left(x\right)=\pi_{n}^{\infty}i^{\left[n\right]}b^{\left[n\right]}\left(x\right)=hb^{\left[n\right]}\left(x\right)\end{eqnarray*}
so that $b_{W}^{\left[n\right]}P_{h}^{\left[n\right]}=hb^{\left[n\right]}$. \end{proof}

What follows is the main result of the paper.
 
\begin{thm}
\label{teo: generated}Every primitively generated braided bialgebra
is isomorphic as a braided bialgebra to the universal enveloping algebra
of its infinitesimal braided Lie algebra. \end{thm}
\begin{proof}
Let $A$ be a primitively generated braided bialgebra with infinitesimal
braided Lie algebra $\left(P,c_{P},b_{P}\right)$ and set $U:=U\left(P,c_{P},b_{P}\right).$
By Theorem \ref{teo: univ U}, the identity map of $P$ can be lifted
to a morphism of braided bialgebras $\varphi:U\rightarrow A.$ Since
$P$ generates $A$ as a $K$-algebra, $\varphi$ is surjective. On
the other hand, by Corollary \ref{cor:P(U)}, $i_{U}:P\rightarrow U$
induces an isomorphism between $P$ and $P\left(U\right).$ Hence
the restriction of $\varphi$ to $P\left(U\right)$ is injective.
By \cite[Lemma 5.3.3]{Montgomery} $\varphi$ is injective too.\end{proof}

\begin{example}\label{ex:classic}
Assume $\mathrm{char}K=0$. Let $V$ be a vector space. Regard it
as a braided vector space through the canonical flip map $c:V\otimes V\rightarrow V\otimes V:x\otimes y\mapsto y\otimes x$.
Then the braiding $c$ is of Hecke type with regular mark $1$ (regularity,
when the mark is $1$, just means $\mathrm{char}K=0$). Thus, by \cite[Theorem 6.14]{Ardizzoni-Sdeg},
we have that $\left(V,c\right)$ is in $\mathcal{S}$. 

Let $b$ such that $\left(V,c,b\right)$ is a braided Lie algebra.
By Theorem \ref{teo:Ardi}, we have that $U\left(V,c,b\right)=U\left(V,c,b^{\left[0\right]}\right)=\mathbb{U}\left(V,c,\beta\right)$,
where $\beta:E\left(V,c\right)\rightarrow V$ is the restriction of
$b^{\left[0\right]}$ to $E\left(V,c\right)$. By Lemma \ref{lem:UnivRank},
we have that, for all $t\geq1$, $b^{\left[t\right]}=b^{\left[1\right]}$.
Moreover $b^{\left[1\right]}$ is bijective and it is the inverse
of $i^{\left[1\right]}:V\rightarrow P^{\left[1\right]}$. By \cite[Theorem 6.3]{Ardizzoni-Universal}\[
\mathbb{U}\left(V,c,\beta\right)=\frac{T\left(V,c\right)}{\left(x\otimes y-c\left(x\otimes y\right)-\left[x,y\right]\mid x,y\in V\right)}=\frac{T\left(V\right)}{\left(x\otimes y-y\otimes x-\left[x,y\right]\mid x,y\in V\right)}\]
where $\left[-,-\right]:V\times V\rightarrow V$ is defined by $\left[x,y\right]:=\beta\left(x\otimes y-c\left(x\otimes y\right)\right)=b^{\left[0\right]}\left(x\otimes y-y\otimes x\right)$
for all $x,y\in V$. It is remarkable that $\left(V,\left[-,-\right]\right)$
is an ordinary Lie algebra (see \cite[Remark 6.4]{Ardizzoni-Universal})
so that $U\left(V,c,b\right)$ is nothing but the ordinary universal
enveloping algebra of $\left(V,\left[-,-\right]\right)$.

Conversely, let $\left[-,-\right]:V\times V\rightarrow V$ be such
that $\left(V,\left[-,-\right]\right)$ is an ordinary Lie algebra.
Let $A$ be the ordinary universal enveloping algebra of $\left(V,\left[-,-\right]\right)$.
This is an ordinary Hopf algebra, see e.g. \cite[Proposition V.2.4, page 97]{Kassel}
or \cite[Proposition 7, page 115]{Bourbaki-LieGrops1-3}. In particular
it is a braided bialgebra with braiding the canonical flip map on
$A$. It is primitively generated as, by construction it is generated
by the image of $V$ in $A$. By Theorem \ref{teo: generated}, $A\cong U\left(P,c_{P},b_{P}\right)$
where $\left(P,c_{P},b_{P}\right)$ is the infinitesimal braided Lie
algebra of $A$. Now, since $\mathrm{char}K=0$, we have that the
canonical map $\sigma:V\rightarrow P=P\left(A\right)$ is bijective
\cite[Corollary at page 117]{Bourbaki-LieGrops1-3}. Since $c_{P}$
is the restriction of the braiding of $A$, then $c_{P}$ is the canonical
flip map $c$ on $V$. Hence we can identify $\left(P,c_{P}\right)$
with $\left(V,c\right)$ and $\sigma$ with the identity of $V$.
Hence we get a braided bialgebra $\left(V,c,b\right)$ where $b=b_{P}$.
By construction the first part we have that, for all $t\geq1$, $b^{\left[t\right]}=b^{\left[1\right]}$.
Moreover $b^{\left[1\right]}$ is the inverse of $\sigma$ so that
$b^{\left[1\right]}=Id_{V}$. On the other hand, one has $b^{\left[0\right]}\left(z\right)=\varphi^{\left[0\right]}\left(z\right)$
for all $P\left(T\right)$ where $T:=T\left(V,c\right)$ and $\varphi^{\left[0\right]}:U^{\left[0\right]}\rightarrow A$
is defined as in the proof of Proposition \ref{pro: alg is Lie}.
Explicitly, if $z=\sum_{0\leq i\leq n}z_{i}$ for some $z_{i}\in V^{\otimes i}$,
then $b^{\left[0\right]}\left(z\right)=\sum_{0\leq i\leq n}\varphi^{\left[0\right]}\left(z_{i}\right)=\sum_{0\leq i\leq n}m_{A}^{i-1}\left(z_{i}\right)$
where $m_{A}^{i-1}:A^{\otimes i}\rightarrow A:x_{1}\otimes\cdots\otimes x_{i}\mapsto x_{1}\cdots x_{i}$
denotes the iterated multiplication of $A$.\end{example}

We are now able to recover the classical Milnor-Moore theorem.

\begin{cor}(cf. \cite[Theorem 5.18]{Milnor-Moore})\label{coro:MM} Let $A$ be a cocommutative ordinary bialgebra in characteristic zero. Then $A$ is isomorphic as a bialgebra to the universal enveloping algebra of $P(A)$.
\end{cor}

\begin{proof}
 In vie of Example \ref{ex: cocommutative}, $A$ is primitively generated. By Theorem \ref{teo: generated}, $A$ is isomorphic as a braided bialgebra to the universal enveloping algebra
of its infinitesimal braided Lie algebra $(P,c_P,b_P)$. By Example \ref{ex:classic}, $U(P,c_P,b_P)$ is the ordinary universal enveloping algebra of the ordinary Lie algebra $(P,[-,-])$, where $\left[x,y\right]:=b_P^{\left[0\right]}\left(x\otimes y-y\otimes x\right)=x\cdot_A y-y\cdot_A x$, for all $x,y\in P$.
\end{proof}

\begin{claim}
\label{cla:idealTower}Let $A$ be a braided bialgebra with infinitesimal
braided Lie algebra $\left(V,c,b\right)$. Let $h=\mathrm{Id}_{V}.$
and let $U:=U\left(V,c,b\right)$. Let $h^{\left[n\right]}:U^{\left[n\right]}\rightarrow A$
be defined as in the proof of Theorem \ref{teo: univ U}. Mimic the
construction of \cite[Section 3]{Kharchenko-Acombinatorial} (see
also \cite[Section 3]{Masuoka- pre Nichols}) as follows. Let $I:=\mathrm{Ker}\left(h^{\left[0\right]}\right)$.
Set $I_{0}:=0$, $I_{1}:=\left(P\left[I\right]\right)$, $I_{2}/I_{1}:=\left(P\left[I/I_{1}\right]\right)$,
$I_{3}/I_{2}:=\left(P\left[I/I_{2}\right]\right)$ and so on. Continuing
this way one gets an increasing chain of ideals of $U^{\left[0\right]}=T\left(V,c\right)$
\[
0=I_{0}\leq I_{1}\leq I_{2}\leq\cdots\leq I_{n}\leq\cdots.\]
\end{claim}
\begin{prop}
Keep the assumptions and notations of \ref{cla:idealTower}. Then,
for all $n\in\mathbb{N}$ \[
\mathrm{Ker}\left(h^{\left[n\right]}\right)=\pi_{0}^{n}\left(I\right)\quad\text{and}\quad\mathrm{Ker}\left(\pi_{0}^{n}\right)=I_{n},\]
where $\pi_{0}^{n}:T\left(V,c\right)=U^{\left[0\right]}\rightarrow U^{\left[n\right]}$
is the canonical projection. Moreover \[
\mathrm{Ker}\left(\pi_{0}^{\infty}\right)=\bigcup_{n\in\mathbb{N}}I_{n}.\]
In particular $I=\bigcup_{n\in\mathbb{N}}I_{n}$ whenever $A$ is primitively
generated.\end{prop}
\begin{proof}
We proceed by induction on $n\in\mathbb{N}$. For $n=0$, there is
nothing to prove. Let $n>1$ and assume the formulas true for $n-1$.
Then, since $U^{\left[n\right]}=U^{\left[n-1\right]}/\mathrm{Ker}\left(\pi_{n-1}^{n}\right)$,
we get\[
\mathrm{Ker}\left(h^{\left[n\right]}\right)=\frac{\mathrm{Ker}\left(h^{\left[n-1\right]}\right)}{\mathrm{Ker}\left(\pi_{n-1}^{n}\right)}=\frac{\pi_{0}^{n-1}\left(I\right)}{\mathrm{Ker}\left(\pi_{n-1}^{n}\right)}=\pi_{n-1}^{n}\pi_{0}^{n-1}\left(I\right)=\pi_{0}^{n}\left(I\right).\]
Let us check that $I_{n}=\mathrm{Ker}\left(\pi_{0}^{n}\right)$. We
have\[
\mathrm{P\left[\pi_{0}^{n-1}\left(I\right)\right]=P\left(\mathrm{Ker}\left(h^{\left[n-1\right]}\right)\right)=P^{\left[n-1\right]}\cap\mathrm{Ker}\left(h^{\left[n-1\right]}\right)=\mathrm{Ker}\left(h_{\mid P^{\left[n-1\right]}}^{\left[n-1\right]}\right)=Ker\left[b^{\left[n-1\right]}\right]}.\]
Thus, by Corollary \ref{coro: ker b}, \begin{equation}
\mathrm{Ker}\left(\pi_{n-1}^{n}\right)=\left(Ker\left[b^{\left[n-1\right]}\right]\right)=\left(P\left[\pi_{0}^{n-1}\left(I\right)\right]\right).\label{eq:formKhar}\end{equation}
Now, since $I_{n-1}=\mathrm{Ker}\left(\pi_{0}^{n-1}\right)$, then
the map \[
\alpha_{n-1}:U^{\left[0\right]}/I_{n-1}\rightarrow U^{\left[n-1\right]}:z+I_{n-1}\mapsto\pi_{0}^{n-1}\left(z\right)\]
 is a well defined braided bialgebra isomorphism. Thus \[
I_{n}=\mathrm{Ker}\left(\pi_{0}^{n}\right)\Leftrightarrow\frac{I_{n}}{I_{n-1}}=\frac{\mathrm{Ker}\left(\pi_{0}^{n}\right)}{I_{n-1}}\Leftrightarrow\alpha_{n-1}\left[\frac{I_{n}}{I_{n-1}}\right]=\alpha_{n-1}\left[\frac{\mathrm{Ker}\left(\pi_{0}^{n}\right)}{I_{n-1}}\right].\]
Since $\alpha_{n-1}$ is a braided bialgebra isomorphism, we get \begin{eqnarray*}
\alpha_{n-1}\left[\frac{I_{n}}{I_{n-1}}\right] & = & \alpha_{n-1}\left[\left(P\left[\frac{I}{I_{n-1}}\right]\right)\right]=\left(\alpha_{n-1}\left[P\left[\frac{I}{I_{n-1}}\right]\right]\right)\\
 & = & \left(P\left[\alpha_{n-1}\left[\frac{I}{I_{n-1}}\right]\right]\right)=\left(P\left[\pi_{0}^{n-1}\left[I\right]\right]\right)\end{eqnarray*}
and\begin{eqnarray*}
\alpha_{n-1}\left[\frac{\mathrm{Ker}\left(\pi_{0}^{n}\right)}{I_{n-1}}\right] & = & \pi_{0}^{n-1}\left[\mathrm{Ker}\left(\pi_{0}^{n}\right)\right]=\pi_{0}^{n-1}\left[\mathrm{Ker}\left(\pi_{n-1}^{n}\pi_{0}^{n-1}\right)\right]\\
 & = & \pi_{0}^{n-1}\left(\pi_{0}^{n-1}\right)^{-1}\left[\mathrm{Ker}\left(\pi_{n-1}^{n}\right)\right]=\mathrm{Ker}\left(\pi_{n-1}^{n}\right)\\
 & \overset{\eqref{eq:formKhar}}{=} & \left(P\left[\pi_{0}^{n-1}\left(I\right)\right]\right)\end{eqnarray*}
Hence $I_{n}=\mathrm{Ker}\left(\pi_{0}^{n}\right)$. Now, we have
\[
\mathrm{Ker}\left(\pi_{0}^{\infty}\right)=\mathrm{Ker}\left(\underrightarrow{\lim}\pi_{0}^{n}\right)=\sum_{n\in\mathbb{N}}\mathrm{Ker}\left(\pi_{0}^{n}\right)=\sum_{n\in\mathbb{N}}I_{n}=\bigcup_{n\in\mathbb{N}}I_{n}.\]
By Theorem \ref{teo: generated}, if $A$ is primitively generated
then there is an isomorphism $h^{\left[\infty\right]}:U\rightarrow A$
such that $h^{\left[\infty\right]}\circ\pi_{0}^{\infty}=h^{\left[0\right]}$.
Thus $\mathrm{Ker}\left(h^{\left[0\right]}\right)=\mathrm{Ker}\left(\pi_{0}^{\infty}\right)$
or, equivalently, one has $I=\bigcup_{n\in\mathbb{N}}I_{n}$.
\end{proof}

\section{The trivial bracket}\label{sec:trivbracket}

In this section we associate a bracket to any braided vector space. We give a condition guaranteeing that the bracket of a given braided Lie algebra is trivial. Finally we check this condition holds on a specific example. 

\begin{example}[\textbf{The triavial bracket.}]
\label{ex: symmetric algebra} We now
give a construction of a trivial braided Lie algebra structure for
any braided vector space $\left(V,c\right)$. Let $\left(U^{\left[0\right]},i^{\left[0\right]}\right)$
be defined as in Definition \ref{def: main construction}. Let $b_{tr}^{\left[0\right]}:P^{\left[0\right]}=P\left(U^{\left[0\right]}\right)\rightarrow V$
be the restriction of the canonical projection $U^{\left[0\right]}\rightarrow V$.
Since this projection is a morphism of braided vector spaces, it is
clear that $b_{tr}^{\left[0\right]}$ fulfills \eqref{form:bracket1} and \eqref{form:bracket2} for $n=0$. As in Definition
\ref{def: main construction} we can build $\left(U^{\left[1\right]},i^{\left[1\right]}\right)$.
Since $b_{tr}^{\left[0\right]}i^{\left[0\right]}=\mathrm{Id}_{V}$ we get
\[
U^{\left[1\right]}=\frac{U^{\left[0\right]}}{\left(\left[\mathrm{Id}-i^{\left[0\right]}b_{tr}^{\left[0\right]}\right]\left[P^{\left[0\right]}\right]\right)}=\frac{U^{\left[0\right]}}{\left(\mathrm{Ker}\left[b_{tr}^{\left[0\right]}\right]\right)}.\]
Note that $\mathrm{Ker}\left[b_{tr}^{\left[0\right]}\right]$ is the
space spanned by primitive elements in $T\left(V,c\right)$ of homogeneous
degree at lest $2$. Now $U^{\left[1\right]}$ is a graded braided
bialgebra. The component of degree one is $V^{\left[1\right]}=\mathrm{Im}\left(i^{\left[1\right]}\right)\cong V$
so that we can consider the canonical projection $U^{\left[1\right]}\rightarrow V$.
Moreover, its restriction $b_{tr}^{\left[1\right]}:P^{\left[1\right]}=P\left(U^{\left[1\right]}\right)\rightarrow V$ fulfills \eqref{form:bracket1} and \eqref{form:bracket2} for $n=1$.
Proceeding this way one gets a bracket $b_{tr}=\left(b_{tr}^{\left[n\right]}\right)_{n\in\mathbb{N}}$
on the braided vector space $\left(V,c\right)$ that will be called
the \textbf{trivial bracket} of $\left(V,c\right)$. By Proposition \ref{pro: i_U split},
$\left(V,c,b_{tr}\right)$ is a braided Lie algebra. The universal enveloping
algebra of rank $n$ in this case will be also denoted by \[
S^{\left[n\right]}:=S^{\left[n\right]}\left(V,c\right):=U\left(V,c,\left(b_{tr}^{\left[i\right]}\right)_{0\leq i\leq n-1}\right)\]
and it will be called the \textbf{symmetric algebra of rank $n$ of
}$\left(V,c\right)$. It is the braided bialgebra $S^{\left[n\right]}\left(B\right)$
introduced in \cite[Definition 3.11]{Ardizzoni-Sdeg} in the case
$B=U^{\left[0\right]}$. Thus, since the tensor algebra $U^{\left[0\right]}$
is strongly $\mathbb{N}$-graded as an algebra, by applying \cite[Corollary A.11]{Ardizzoni-Sdeg},
we obtain that \[
U\left(V,c,b_{r}\right)=\underrightarrow{\lim}U^{\left[n\right]}\left(V,c,\left(b_{tr}^{\left[i\right]}\right)_{0\leq i\leq n-1}\right)=\underrightarrow{\lim}S^{\left[n\right]}=S^{\left[\infty\right]}\]
 identifies with the Nichols algebra $\mathcal{B}\left(V,c\right)$.\end{example}
\begin{prop}
\label{pro:vanish}Let $\left(V,c,b\right)$ be a braided Lie algebra.
Let $p_{0}^{\infty}:T\left(V,c\right)\rightarrow\mathcal{B}\left(V,c\right)$
and $\pi_{0}^{\infty}:T\left(V,c\right)\rightarrow U\left(V,c,b\right)$
be the canonical projections. Suppose the ideal $\mathrm{Ker}\left(p_{0}^{\infty}\right)$
is generated by a set $W$. The following are equivalent.
\begin{itemize}
\item[$(i)$] $\pi_{0}^{\infty}$ vanishes
on $W.$
\item[$(ii)$] $b$ is the trivial bracket of $\left(V,c\right).$
\end{itemize}
\end{prop}
\begin{proof}
$\left(i\right)\Rightarrow\left(ii\right)$ Since $\mathrm{Ker}\left(p_{0}^{\infty}\right)=\left(W\right)$
and $\pi_{0}^{\infty}$ vanishes on $W$, there exists a unique algebra
homomorphism $\phi:\mathcal{B}\left(V,c\right)\rightarrow U(V,c,b)$
such that $\phi\circ p_{0}^{\infty}=\pi_{0}^{\infty}$. From the surjectivity
of $p_{0}^{\infty}$ one deduce that $\phi$ is indeed a braided bialgebra
homomorphism. Moreover $\phi$ is surjective as $\pi_{0}^{\infty}$
is. Note that \[
\phi_{\mid V}=\left(\phi\circ p_{0}^{\infty}\right)_{\mid V}=\left(\pi_{0}^{\infty}\right)_{\mid V}=i_{U}\]
which is injective. Since $P\left(\mathcal{B}\left(V,c\right)\right)\cong V$,
by \cite[Lemma 5.3.3]{Montgomery}, $\phi$ is injective whence bijective.

This isomorphism yields an isomorphism of braided brackets between
the infinitesimal part of $\mathcal{B}\left(V,c\right)=U\left(V,c,b_{r}\right)$
and the infinitesimal part of $U(V,c,b)$. In view of Theorem
\ref{teo:infpartU}, we get that there is an isomorphism of braided
brackets between $h:\left(V,c,b_{r}\right)\rightarrow(V,c,b)$. By
definition of $\phi$ it is also clear that $h=\mathrm{Id}_{V}$.
By \eqref{form:morphbrackets}, for each $n\in\mathbb{N}$, we get
$b_{r}^{\left[n\right]}=hb_{r}^{\left[n\right]}=b^{\left[n\right]}P_{h}^{\left[n\right]}=b^{\left[n\right]}$
so that $b=b_{r}.$ %

$\left(ii\right)\Rightarrow\left(i\right)$ Since $b$ is the trivial
bracket of $\left(V,c\right),$ by Example \ref{ex: symmetric algebra},
we have that $U\left(V,c,b\right)=\mathcal{B}\left(V,c\right)$ and
$\pi_{0}^{\infty}=p_{0}^{\infty}$. 
\end{proof}

We now give an example of a braided vector space $\left(V,c\right)$
which becomes a braided Lie algebra only in a trivial way. It is remarkable
that $\left(V,c\right)$ does not belong to the class $\mathcal{S}$
of Remark \ref{rem:Ardi}. Indeed this gives a counterexample to Theorem
\ref{teo:Ardi}. In fact \[
U\left(V,c,b_{tr}\right)=\mathcal{B}\left(V,c\right)\neq S^{\left[1\right]}\left(V,c\right)=U\left(V,c,b_{tr}^{\left[0\right]}\right).\]
 
\begin{example}\label{ex:combRank2}
Let $K$ be a fixed field with $\mathrm{char}\left(K\right)\neq2.$
At the end of \cite{Kharchenko-SkewPrim}, an example of a two dimensional
braided vector space $\left(V,c\right)$ of combinatorial rank $2$
is given. This essentially means that $S^{\left[0\right]}\left(V,c\right)\neq S^{\left[1\right]}\left(V,c\right)\neq S^{\left[2\right]}\left(V,c\right)=\mathcal{B}\left(V,c\right)$.
The braiding $c$ is of diagonal type of the form\[
c\left(x_{i}\otimes x_{j}\right)=q_{i,j}x_{j}\otimes x_{i},1\leq i,j\leq2.\]
where $x_{1},x_{2}$ is a basis of $V$ over $K$, $q_{1,2}=1\neq-1$
and $q_{i,j}=-1$ for all $\left(i,j\right)\neq\left(1,2\right).$
One has that \[
\mathcal{B}\left(V,c\right)=\frac{T\left(V,c\right)}{\left(x_{1}^{2},\quad x_{2}^{2},\quad x_{1}x_{2}x_{1}x_{2}+x_{2}x_{1}x_{2}x_{1}\right)}.\]
Thus $\mathcal{B}\left(V,c\right)$ has basis \[
1_{K},x_{1},x_{2},x_{2}x_{1},x_{1}x_{2},x_{1}x_{2}x_{1},x_{2}x_{1}x_{2},x_{2}x_{1}x_{2}x_{1}.\]
Note that the endomorphism $c:V\otimes V\rightarrow V\otimes V$ has
minimal polynomial $\left(X+1\right)\left(X^{2}+1\right)$. Our aim
here is to determine all primitively generated connected braided bialgebras
$A$ with infinitesimal part $\left(V,c\right)$. Let $A$ be such
a braided bialgebra. By Theorem \ref{teo: generated}, $A$ is isomorphic
to the enveloping algebra $U\left(P,c_{P},b_{P}\right)$ associated
to its infinitesimal braided Lie algebra $\left(P,c_{P},b_{P}\right)$.
Therefore we have to investigate the possible braided Lie algebra
structures on $\left(V,c\right)$. 

Let $b$ be a bracket on $\left(V,c\right)$ such that $\left(V,c,b\right)$
is a braided Lie algebra. Set $T:=T\left(V,c\right).$ Clearly $x_{1}^{2},x_{2}^{2}\in P^{\left[0\right]}=P\left(T\right).$ 

Suppose there is $t\in\left\{ 1,2\right\} $ such that $b^{\left[0\right]}\left(x_{t}^{2}\right)\ne0$.
Note that $x_{t}^{2}$ is a homogeneous quantum operation in the sense
of \cite[Definition 9.1]{Ardizzoni-Universal}. By \cite[Corollary 9.6]{Ardizzoni-Universal},
there is $a\in\left\{ 1,2\right\} $ such that the following relation
holds \begin{eqnarray*}
q_{a,j} & = & \prod_{l=1,2}q_{l,j}^{d_{l}}=q_{t,j}^{2}=1,\text{ for all }j=1,2.\end{eqnarray*}
In particular $-1=q_{a,1}=1$ which contradicts the hypothesis $\mathrm{char}\left(K\right)\neq2.$
Hence we proved that \[
b^{\left[0\right]}\left(x_{t}^{2}\right)=0\text{ for }t=1,2.\]
Therefore\begin{equation}
\pi_{0}^{1}\left(x_{t}^{2}\right)=\pi_{0}^{1}\left[i^{\left[0\right]}b^{\left[0\right]}\left(x_{t}^{2}\right)\right]=0,\text{ for }t=1,2.\label{eq:esK}\end{equation}
Using \eqref{eq:esK}, one easily checks that $w:=\pi_{0}^{1}\left(x_{1}x_{2}x_{1}x_{2}+x_{2}x_{1}x_{2}x_{1}\right)\in P^{\left[1\right]}=P\left(U^{\left[1\right]}\right)$.

Let us check that $b^{\left[1\right]}\left[w\right]=0$. Set $b^{\left[1\right]}\left(w\right)=b_{1}^{\left[1\right]}x_{1}+b_{2}^{\left[1\right]}x_{2}$.
By (\ref{form:bracket1}), for $1\leq j\leq n$ we have that\[
c_{V^{\left[1\right]}}\left(i^{\left[1\right]}b^{\left[1\right]}\otimes V^{\left[1\right]}\right)=\left(V^{\left[1\right]}\otimes i^{\left[1\right]}b^{\left[1\right]}\right)c_{P^{\left[1\right]},V^{\left[1\right]}},\]
Now\begin{eqnarray*}
c_{V^{\left[1\right]}}\left(i^{\left[1\right]}b^{\left[1\right]}\otimes V^{\left[1\right]}\right)\left(w\otimes\pi_{0}^{1}\left(x_{j}\right)\right) & = & c_{V^{\left[1\right]}}\left(i^{\left[1\right]}b^{\left[1\right]}\left(w\right)\otimes\pi_{0}^{1}\left(x_{j}\right)\right)\\
 & = & b_{1}^{\left[1\right]}c_{V^{\left[1\right]}}\left(\pi_{0}^{1}\left(x_{1}\right)\otimes\pi_{0}^{1}\left(x_{j}\right)\right)+b_{2}^{\left[1\right]}c_{V^{\left[1\right]}}\left(\pi_{0}^{1}\left(x_{2}\right)\otimes\pi_{0}^{1}\left(x_{j}\right)\right)\\
 & = & \left(\pi_{0}^{1}\otimes\pi_{0}^{1}\right)\left[b_{1}^{\left[1\right]}c\left(x_{1}\otimes x_{j}\right)+b_{2}^{\left[1\right]}c\left(x_{2}\otimes x_{j}\right)\right]\\
 & = & q_{1,j}b_{1}^{\left[1\right]}\pi_{0}^{1}\left(x_{j}\right)\otimes\pi_{0}^{1}\left(x_{1}\right)+q_{2,j}b_{2}^{\left[1\right]}\pi_{0}^{1}\left(x_{j}\right)\otimes\pi_{0}^{1}\left(x_{2}\right)\\
 & = & q_{1,j}b_{1}^{\left[1\right]}\pi_{0}^{1}\left(x_{j}\right)\otimes\pi_{0}^{1}\left(x_{1}\right)-b_{2}^{\left[1\right]}\pi_{0}^{1}\left(x_{j}\right)\otimes\pi_{0}^{1}\left(x_{2}\right)\end{eqnarray*}
and\begin{eqnarray*}
\left(V^{\left[1\right]}\otimes i^{\left[1\right]}b^{\left[1\right]}\right)c_{P^{\left[1\right]},V^{\left[1\right]}}\left(w\otimes\pi_{0}^{1}\left(x_{j}\right)\right) & = & q_{1,j}^{2}q_{2,j}^{2}\left(V^{\left[1\right]}\otimes i^{\left[1\right]}b^{\left[1\right]}\right)\left(\pi_{0}^{1}\left(x_{j}\right)\otimes w\right)\\
 & = & b_{1}^{\left[1\right]}\pi_{0}^{1}\left(x_{j}\right)\otimes\pi_{0}^{1}\left(x_{1}\right)+b_{2}^{\left[1\right]}\pi_{0}^{1}\left(x_{j}\right)\otimes\pi_{0}^{1}\left(x_{2}\right).\end{eqnarray*}
Since $\pi_{0}^{1}\left(x_{i}\right)=i^{\left[0\right]}\left(x_{j}\right)$
and $i^{\left[0\right]}$ is injective, we deduce that\[
q_{1,j}b_{1}^{\left[1\right]}x_{j}\otimes x_{1}-b_{2}^{\left[1\right]}x_{j}\otimes x_{2}=b_{1}^{\left[1\right]}x_{j}\otimes x_{1}+b_{2}^{\left[1\right]}x_{j}\otimes x_{2}\]
whence, for $j=1$ we get \[
2b_{1}^{\left[1\right]}=0=2b_{2}^{\left[1\right]}\]
whence $b^{\left[1\right]}\left[w\right]=0$ as $\mathrm{char}\left(K\right)\neq2.$
Now\begin{eqnarray*}
\pi_{0}^{\infty}\left(x_{1}x_{2}x_{1}x_{2}+x_{2}x_{1}x_{2}x_{1}\right) & = & \pi_{1}^{\infty}\pi_{0}^{1}\left(x_{1}x_{2}x_{1}x_{2}+x_{2}x_{1}x_{2}x_{1}\right)\\
 & = & \pi_{1}^{\infty}\left(w\right)\\
 & = & \pi_{2}^{\infty}\pi_{1}^{2}\left(w\right)\\
 & \overset{w\in P^{\left[1\right]}}{=} & \pi_{2}^{\infty}\pi_{1}^{2}i^{\left[1\right]}b^{\left[1\right]}\left(w\right)=0.\end{eqnarray*}
Since we already proved that $\pi_{0}^{1}\left(x_{t}^{2}\right)=0$,
we also have $\pi_{0}^{\infty}\left(x_{t}^{2}\right)=0$ for $t=1,2$.
By Proposition \ref{pro:vanish}, $b$ is the trivial bracket of $\left(V,c\right)$.
As a consequence, we get $U\left(V,c,b\right)=\mathcal{B}\left(V,c\right)$.
In view of Theorem \ref{teo: generated}, $\mathcal{B}\left(V,c\right)$
is then the unique primitively generated connected braided bialgebra
$A$ with infinitesimal part $\left(V,c\right)$. 
\end{example}

\end{document}